\documentclass{article}
\usepackage[utf8]{inputenc}
\usepackage{mathrsfs}
\usepackage{amssymb}
\usepackage{amsmath}
\usepackage{bm}
\usepackage{cite}
\usepackage{authblk}
\usepackage{algorithm}
\usepackage{algpseudocode}
\usepackage{tikz,pgfplots}
\pgfplotsset{compat=1.18}
\usepackage{amsfonts}
\usepackage{subfigure}
\usepackage{graphicx}       
\usepackage{braket}
\usepackage{booktabs}
\usepackage{multirow}
\usepackage[margin=1in]{geometry}

\usepackage[normalem]{ulem}
\usepackage[colorlinks,linkcolor=blue,anchorcolor=blue,citecolor=red]{hyperref}
\usepackage{enumerate}
\usepackage{cases}
\usepackage{mathrsfs}
\usepackage{calligra}
\usepackage{makecell}
\newtheorem{remark}{Remark}

\newcommand\bx{\boldsymbol{x}}

\newcommand\bu{\boldsymbol{u}}
\newcommand\bv{{\boldsymbol{v}}}

\newcommand\bc{\boldsymbol{c}}

\newcommand\bbR{\mathbb{R}}

\newcommand\bfS{\mathbf{S}}

\newcommand\bq{\boldsymbol{q}}

\newcommand\Kn{{\rm Kn}}

\newcommand\mQ{\mathcal{Q}}
\newcommand\mM{\mathcal{M}}

\newcommand\mP{\mathcal{P}}

\newcommand\pd[2]{\frac{\partial {#1}}{\partial {#2}}}

\newcommand\mL{\mathcal{L}}
\newcommand{\mO}{\mathcal{O}}

\newcommand{\tbg}{\tilde{\bm{g}}}

\newcommand{\dd}{\; \mathrm{d}}

\newcommand\hatn[1]{{#1}^{(n)}}
\graphicspath{{media/}}     

\newcommand{\tf}[1]{\left(#1\right)_{\mathrm{stf}}}
\numberwithin{equation}{section}
\newcommand{\NA}{\multicolumn{1}{c}{-}}

\title{Accelerated iterative method for solving
the steady-state Boltzmann equation}
\date{}
\author{
  ~~Pei Zhang\footnote{Beijing Computational Science Research Center,
  Beijing, China, email:
  \texttt{zhangpei@csrc.ac.cn}.},
  ~~Zhenning Cai\footnote{Department of Mathematics, National University of Singapore, Singapore, 119076, email:
  \texttt{matcz@nus.edu.sg}.},
  ~~Yanli Wang\footnote{Beijing Computational Science Research Center,
  Beijing, China, email:
  \texttt{ylwang@csrc.ac.cn}.}
}
\begin{document}
\maketitle

\begin{abstract}
The efficient simulation of steady-state rarefied gas flows remains a significant computational challenge due to the high dimensionality of the collision integral and the severe numerical stiffness in the near-continuum regime. In this work, we propose a modified Newton method equipped with a macroscopic synthetic system (Newton-MS) for the steady-state Boltzmann equation with the quadratic collision operator. In Newton-MS, the modified Newton iteration is utilized as the outer nonlinear solver, while each Newton correction equation is solved by an inner source iteration, where the linearized collision operator is utilized to approximate the quadratic collision model, and it is reduced into a linear iteration. Moreover, a macroscopic synthetic system based on Chapman-Enskog closure is derived to accelerate the convergence of the linear inner iteration in the continuum limit. 
Besides, the fully discrete macroscopic synthetic system is deduced under the framework of the discontinuous Galerkin method to reduce computational cost compared to directly discretizing the continuous macroscopic synthetic system. Several numerical examples, including the 1D Fourier, Couette flow problem, and the 2D cavity flow and thermal-driven cavity flow, are studied to validate the high efficiency of Newton-MS.

    
\end{abstract}

\section{Introduction}
\label{intro}

Rarefied gas dynamics constitutes a pivotal branch of fluid mechanics, essential for understanding multiscale flow phenomena in high-altitude aerothermodynamics, micro-electromechanical systems (MEMS), and specialized gas transport processes. The flow regimes are characterized by the Knudsen number ($\Kn$), defined as the ratio of the molecular mean free path to a characteristic macroscopic length scale. Although the Euler or Navier-Stokes equations suffice for the continuum regime ($\Kn \ll 1$), they become physically inadequate as the rarefaction effects increase. Consequently, the Boltzmann equation, which governs the evolution of the single-particle probability distribution function in phase space, serves as the fundamental framework that is valid across all flow regimes. However, the high dimensionality of phase space and the complex quadratic collision term pose substantial challenges for efficient numerical simulations.

Traditional numerical methods for solving the Boltzmann equation are generally categorized into stochastic methods and deterministic methods. Direct Monte Carlo simulation (DSMC) \cite{bird1994molecular, oran1998direct, ganjaei2009new} is one of the most popular stochastic methods that is effective in simulating high-speed rarefied gas flows. However, it suffers from computational inefficiency and inherent statistical noise in low-speed regimes. The Unified Gas-Kinetic Wave-Particle (UGKWP) method couples deterministic waves for continuum flows with stochastic particles for kinetic non-equilibrium, covering the full Knudsen number regime. Unlike DSMC, which is costly and time-step limited in near-continuum flows, UGKWP adaptively weights waves and particles by local Knudsen number, recovering a fluid solver in continuum and a particle method in rarefied regimes. This hybrid reduces noise and cost while preserving kinetic accuracy for multiscale gas dynamics \cite{liu2020unifiedI}.

Conversely, deterministic methods are often preferred for continuous or low-speed flows where the statistical fluctuations must be minimized. The discrete velocity method (DVM) \cite{broadwell1964study, buet1996discrete, liu2020unified} is a classical deterministic method, which discretizes the distribution function at a series of microscopic velocity points.  The unified gas-kinetic scheme (UGKS) can eliminate restrictions on cell size and time step by simultaneously handling free flow and collision of gas molecules \cite{liu2016unified, liu2020unifiedI}. However, since information exchange depends on the evolution of the velocity distribution function, UGKS still requires a large number of iterations to obtain the steady-state solution of the near-continuum flow \cite{chen2015comparative, wang2018comparative}. Furthermore, spectral methods, including the Fourier spectral method \cite{pareschi1996fourier, pareschi2000numerical, gamba2017fast}, the Hermite spectral method \cite{hu2020numerical, kitzler2019polynomial}, the Burnett spectral method \cite{hu2020burnett, gamba2018galerkin}, and the mapped Chebyshev spectral method \cite{hu2022fast, hu2020petrov} utilize global basis functions to approximate the distribution function, leading to a higher order of convergence. Alternatively, the moment method \cite{grad1949kinetic} originally proposed by Grad simplifies the kinetic description by approximating the distribution function through a truncated series of macroscopic moments.

For the steady-state problems, which are of great engineering interest, the efficiency of the numerical solver is critical. The conventional source iteration (SI) is efficient in highly rarefied regimes, but exhibits prohibitively slow convergence as the flow approaches the continuum limit ($\Kn \to 0$) \cite{valougeorgis2003acceleration}. In this diffusive limit,  the ``converged" solution is often contaminated by numerical errors arising from velocity discretization and collision operator approximations. To overcome the slow convergence in the near-continuum regime, Synthetic Iterative Schemes (SIS) have been developed to accelerate the solution of the Boltzmann equation \cite{lihnaropoulos2007formulation, willian2014moment}. The core mechanism of SIS is to couple the microscopic kinetic transport with macroscopic fluid equations, where the macroscopic fluid equations are utilized to guide the evolution of the distribution function. While early SIS versions were often limited to specific flow scenarios or linearized equations, the subsequently proposed General Synthetic Iterative Scheme (GSIS) extended this capability to general rarefied gas flows \cite{su2020can, su2020fast, zhu2021general}. By rigorously deriving the macroscopic equations from the kinetic level, GSIS ensures asymptotic preserving properties and achieves fast convergence across the entire range of Knudsen numbers, typically converging within dozens of iterations.

Alternatively, the Newton method addresses the collision nonlinearity by solving a sequence of Newton equations, each of which is essentially a linearized Boltzmann equation with a residual source term \cite{yin2025fast}. The competitiveness of this method has been substantially enhanced by the recent developments in the fast Fourier spectral method (FFT), where the complexity of the linearized collision operator is reduced to $\mathcal{O}(N^4 \log N)$, with $N$ denoting the number of modes in each velocity direction \cite{yin2025fast}. The computational cost is significantly lower than that of the quadratic collision operator with complexity of $\mathcal{O}(MN^4 \log N)$ \cite{gamba2017fast}, where $M$ is the number of quadrature points on the unit sphere. This new efficient FFT to solve the linear collision operator makes the Newton iteration a highly viable and efficient choice.

However, efficiently evaluating the collision operator does not equate to efficiently solving the linear system. Although the linearized collision operator can be evaluated efficiently, the corresponding source iteration still converges slowly for small $\Kn$, where macroscopic information is propagated inefficiently only through the update of the kinetic equation. This observation motivates a macroscopic synthetic preconditioning strategy, in which a macroscopic system is utilized to accelerate the slow macroscopic components of the Newton correction. However, the macroscopic equations here cannot be directly borrowed from the existing GSIS. The key difference lies in the state around which the equation is linearized. Precisely, nonlinear GSIS is built for the original nonlinear equation \cite{zhu2021general}, and linearized GSIS is derived around a global equilibrium \cite{su2020fast}, which is spatially uniform. For the Newton correction, it is linearized around a local Maxwellian that varies in space and is updated at each outer step. As a result, a different form of the macroscopic synthetic system should be derived specifically.

In this work, a Newton method equipped with a macroscopic synthetic system (Newton-MS) is developed 
for the steady-state Boltzmann equation with the quadratic collision operator. The main idea is retaining the modified Newton iteration for the nonlinear collision term, while constructing a macroscopic synthetic system for the inner correction equation to accelerate the convergence of the slow macroscopic components. As discussed above, the corresponding macroscopic equations differ from existing GSIS and are derived for the Newton correction around a spatially varying local Maxwellian. These continuous equations, however, contain several background-dependent terms induced by the spatial variation of the local Maxwellian, which make direct discretization cumbersome. Rather than discretizing them directly, the matrix form of the macroscopic system at the fully discrete level is constructed, starting from the discrete Newton correction equation. This yields a compact macroscopic linear system that can be assembled once per Newton step and reused in the inner iterations, without explicitly deducing the complicated continuous macroscopic synthetic systems.

The spatial discretization is achieved by the discontinuous Galerkin (DG) method, ensuring high-order accuracy and geometric flexibility for multidimensional problems. Extensive numerical experiments, ranging from the 1D Fourier and Couette flows to 2D cavity flows, are presented to validate the efficiency and accuracy of the proposed Newton-MS method.

The rest of this paper is organized as follows. Sec. \ref{sec:Bolt} briefly reviews the Boltzmann equation and the fast Fourier spectral method for the collision operator. The construction of the Newton-MS, especially the deduction of the continuous macroscopic synthetic system, is presented in Sec. \ref{sec:method}, while Sec. \ref{sec:num_dis} provides the fully discrete form of Newton-MS. Several numerical examples are provided in Sec. \ref{sec:num} with some concluding remarks and the appendix listed in Sec. \ref{sec:con} and App. \ref{sec:constitutive_derivation}.

\section{Preliminaries}
\label{sec:Bolt}
In this section, several basic properties of the steady-state Boltzmann equation are introduced, and the fast Fourier spectral method, which is utilized to approximate the quadratic collision model, is proposed in Sec. \ref{sec:fsm}.
\subsection{The steady-state Boltzmann equation}
\label{sec:ssBolt}
The dynamics of a rarefied gas are characterized by the evolution of the particle distribution function $f(t, \bx, \bv)$ in the phase space, where $t \geqslant 0$ is time, $\bx \in \Omega \subset \mathbb{R}^d,d=1,2,3$ is position, and $\bv \in \mathbb{R}^3$ is the microscopic velocity of gas molecules. Here, the steady-state Boltzmann equation is considered as 
\begin{equation} \label{eq:Boltz}
     \bv \cdot \nabla_{\bx} f(\bx,\bv) = \frac{1}{\Kn}\mathcal{Q}[f, f](\bx,\bv),
\end{equation}
where $\Kn$ is the Knudsen number and $\mathcal{Q}[f, f]$ is the quadratic collision operator, which describes binary collisions between particles and has the following integral form
\begin{equation}
\label{eq:quad_col}
    \mathcal{Q}[f, f](\bx,\bv) = \int_{\mathbb{R}^3} \int_{\mathbb{S}^2} \mathcal{B}(\bv - \bv_*, \boldsymbol{\sigma}) \left[ f(\bx,\bv') f(\bx,\bv'_*) - f(\bx,\bv) f(\bx,\bv_*) \right] \dd \boldsymbol{\sigma} \dd \bv_*,
\end{equation}
where $(\bv, \bv_*)$ and $(\bv', \bv'_*)$ denote the pre- and post-collision velocity pairs, respectively. 
Under the constraints of momentum and energy conservation, the post-collision velocities are parameterized by the unit vector $\boldsymbol{\sigma} \in \mathbb{S}^2$
\begin{equation}
    \bv' = \frac{\bv + \bv_*}{2} + \frac{|\bv - \bv_*|}{2}\boldsymbol{\sigma}, \quad \bv'_* = \frac{\bv + \bv_*}{2} - \frac{|\bv - \bv_*|}{2}\boldsymbol{\sigma}.
\end{equation}
The collision kernel $\mathcal{B}$ reflects the intermolecular potential between gas molecules. A widely adopted simplification is the variable hard-sphere (VHS) model \cite{bird1994molecular}, which assumes $\mathcal{B}$ is independent of the scattering angle $\boldsymbol{\sigma}$
\begin{equation}
\mathcal{B} = C |\bv - \bv_*|^\alpha, \quad \alpha = 2(1 - \omega),
\end{equation}
where $\omega$ represents the viscosity index of the gas and the constant $C$ is determined by the molecular diameter. Due to the complex form of the quadratic collision, several simplified collision models are proposed, such as the linearized collision model
\begin{equation}
\label{eq:linear_col}
    \mathcal{L}[g] = \mathcal{Q}[\mathcal{M},g]+\mathcal{Q}[g,\mathcal{M}],
\end{equation}
where $g$ denotes the distribution acted on by the linearized collision operator. The Maxwellian $\mathcal{M} = \mathcal{M}_{[\rho,\bu,T]}(\bv) $ is given by
\begin{equation}
\label{eq:Maxwell}
    \mathcal{M}_{[\rho,\bu,T]}(\bv) = \frac{\rho}{\sqrt{2\pi T}^{3}}\exp\left(-\frac{|\bv-\boldsymbol{u}|^2}{2T}\right).
\end{equation}
 Here, $\rho$ is the density, $\boldsymbol{u}=(u_1,u_2,u_{3})^T$  is the macroscopic velocity, and $T$ is the temperature. When the reference Maxwellian is associated with a distribution function $f$, the corresponding macroscopic quantities are given by
\begin{equation}
    \label{eq:mac_var}
       \rho = \int_{\bbR^{3}} f(\bx,\bv)\dd\bv,\qquad \boldsymbol{u} =\frac{1}{\rho} \int_{\bbR^{3}} \bv f(\bx,\bv) \dd\bv,\qquad T= \frac{1}{3\rho}\int_{\bbR^{3}} |\bv-\bu|^2 f(\bx,\bv) \dd\bv.
\end{equation}


\subsection{Fast Fourier spectral approximation of the collision operator}
\label{sec:fsm}
Due to the complex form of the quadratic collision model $\mathcal{Q}[f, f]$, numerical methods for it have been extensively studied in recent years. The fast Fourier spectral method (FFS) \cite{pareschi2000numerical, gamba2017fast} is adopted here, due to its high efficiency and easy implementation. As stated in \cite{gamba2017fast}, the time complexity for FFS to approximate the quadratic collision operator is $\mathcal{O}(MN^4\log N)$. For VHS models, the linearized collision operator $\mathcal{L}[g]$ can be approximated with time complexity $\mathcal{O}(N^4\log N)$, according to the recent research work \cite{yin2025fast}. In this section, these algorithms will be briefly reviewed.

For FFS, it truncates the microscopic velocity space to a finite domain $\mathcal{D}_v = [-L, L]^3$. Let $N$ denote the one-dimensional spectral
truncation parameter, so that $2N$ velocity grid points are used in each
direction. The distribution function $f(\bx,\bv)$ is approximated by
\begin{equation}
\label{eq:dis_f}
f(\bv)=\frac{1}{(2L)^3}\sum_{\boldsymbol{k}=-N}^{N}\hat{f}_{\boldsymbol{k}}\,e^{i\pi \boldsymbol{k} \cdot \bv/L},
\qquad
\hat{f}_{\boldsymbol{k}}=\frac{1}{c_{\boldsymbol{k}}}\left(\frac{L}{N}\right)^3\sum_{\boldsymbol{\ell}=-N}^{N-1}
f\!\left(\frac{\boldsymbol{\ell} L}{N}\right)e^{-i\pi \boldsymbol{k}\cdot \boldsymbol{\ell}/N},
\end{equation}
where $\boldsymbol{k} = (k_1, k_2, k_3)$ and
\begin{equation}
\sum_{\boldsymbol{k}=-N}^{N} = \sum_{k_1=-N}^{N} \sum_{k_2=-N}^{N} \sum_{k_3=-N}^{N}, \quad c_{\boldsymbol{k}} = c_{k_1} c_{k_2} c_{k_3}, \quad c_j = \begin{cases} 2, & \text{if } j = \pm N, \\ 1, & \text{otherwise.} \end{cases}    
\end{equation}
Moreover, the discretization of the quadratic collision model $\mathcal{Q}[f, f]$ requires a truncation of the collision kernel by assuming $\mathcal{B}(\bv - \bv_*, \boldsymbol{\sigma}) = 0$ when $|\bv - \bv_*| > R$, where $R > 0$ is a problem-dependent parameter \cite{gamba2017fast}. Then the Fourier coefficients of $\mathcal{Q}$ are calculated as 
\begin{equation}
\label{eq:coe_Q}
\widehat{\mathcal{Q}}_{\boldsymbol{k}} = \sum_{j=1}^J \sum_{m=1}^M C_{j,m,\boldsymbol{k}} \sum_{\boldsymbol{l} = -N}^{N-1} (\alpha^{(1)}_{j,m,\boldsymbol{k}-\boldsymbol{l}} \beta^{(1)}_{j,m,\boldsymbol{l}} + \alpha^{(2)}_{j,m,\boldsymbol{k}-\boldsymbol{l}} \beta^{(2)}_{j,m,\boldsymbol{l}}) \hat{f}_{\boldsymbol{k}-\boldsymbol{l}} \hat{f}_{\boldsymbol{l}} - \sum_{\boldsymbol{l} = -N}^{N-1} \gamma_{\boldsymbol{l}} \hat{f}_{\boldsymbol{k}-\boldsymbol{l}} \hat{f}_{\boldsymbol{l}},
\end{equation}
where all coefficients $C_{j,m,\boldsymbol{k}}$, $\alpha^{(1)}_{j,m,\boldsymbol{k}}$, $\alpha^{(2)}_{j,m,\boldsymbol{k}}$, $\beta^{(1)}_{j,m,\boldsymbol{k}}$, $\beta^{(2)}_{j,m,\boldsymbol{k}}$ and $\gamma_{\boldsymbol{l}}$ can be precomputed, and we refer to \cite{gamba2017fast} for the detailed calculation of these coefficients. Thus, if the fast Fourier transform (FFT) is adopted to compute the convolution, the final computational cost of \eqref{eq:coe_Q} is $\mathcal{O}(JM N^3 \log N)$, where $J$ is the number of quadrature points for a one-dimensional integral over $[0,R]$. Generally speaking, $R$ is chosen to be proportional to $L$, and it is expected that $J = \mathcal{O}(N)$. Thus the computational cost can also be written as $\mathcal{O}(M N^4 \log N)$ \cite{gamba2017fast}.

For the VHS models, the Fourier coefficients of the linearized collision operator $\mathcal{L}[g]$ can be computed more efficiently. The algorithm in \cite{yin2025fast} provides the following form of $\widehat{\mathcal{L}}_{\boldsymbol{k}}$:
\begin{equation}
\label{eq:coe_L}
\widehat{\mathcal{L}}_{\boldsymbol{k}} = \sum_{j=1}^J \alpha_{j,\boldsymbol{k}} \sum_{\boldsymbol{l} = -N}^{N-1} \mathcal{M}_{\boldsymbol{l}}^h \psi_{j,\boldsymbol{l}} e^{-i\pi \boldsymbol{k} \cdot \boldsymbol{l}/N} - \sum_{\boldsymbol{l} = -N}^{N-1} \gamma_{\boldsymbol{l}} (\hat{f}_{\boldsymbol{k}-\boldsymbol{l}} \hat{\mathcal{M}}_{\boldsymbol{l}} + \hat{\mathcal{M}}_{\boldsymbol{k}-\boldsymbol{l}} \hat{f}_{\boldsymbol{l}}),
\end{equation}
where $\hat{\mathcal{M}}_{\boldsymbol{l}}$ is the Fourier coefficient of the local Maxwellian \eqref{eq:Maxwell}, and
{\small
\begin{gather}
\psi_{j,\boldsymbol{l}} = \left(\frac{\pi}{N}\right)^3 \sum_{\boldsymbol{m}} \beta_{j,\boldsymbol{m}} r_{\boldsymbol{l}-\boldsymbol{m}} \mathcal{M}^c_{\boldsymbol{m}}, \qquad
r_{\boldsymbol{l}} = g\left( \frac{\boldsymbol{l}L}{N} \right) \Bigg/\mathcal{M}\left( \frac{\boldsymbol{l}L}{N} \right), \\
\mathcal{M}_{\boldsymbol{m}}^h = \frac{\rho}{(\pi T)^{3/2}} \exp \left( -\frac{1}{T} \left| \frac{\boldsymbol{m} L}{N} - \boldsymbol{u} \right|^2 \right), \qquad
\mathcal{M}_{\boldsymbol{m}}^c = \frac{\rho}{(\pi T)^{3/2}} e^{-\frac{|\boldsymbol{m}|^2 L^2}{N^2 T}}.
\end{gather}
}
The coefficients $\alpha_{j,\boldsymbol{k}}$ and $\beta_{j,\boldsymbol{m}}$ can again be precomputed. With FFT, this leads to a total computational cost $\mathcal{O}(N^4 \log N)$, removing the dependence on $M$, which is the degree of freedom for the integration on the unit sphere compared with the quadratic collision operator.

With the lower computational cost of the linearized collision operator, one can achieve a cheaper model by replacing $\mathcal{Q}[f,f]$ with $\mathcal{L}[g]$, while still maintaining much better accuracy than BGK-type models. Moreover, as briefly tested in \cite{yin2025fast}, the lower computational cost of $\mathcal{L}[g]$ can also help develop efficient algorithms for the steady-state Boltzmann equation \eqref{eq:Boltz}, whose efficiency will be further improved in the following sections, especially in the case of small Knudsen number.

\section{Construction of Newton-macroscopic synthetic system}
\label{sec:method}
Solving the steady-state Boltzmann equation efficiently is hindered by two primary challenges: the prohibitive computational cost of the high-dimensional collision integral, and the slow convergence caused by the stiffness of the transport term in the near-continuum regime. To address these issues, a Newton iteration accelerated by a macroscopic synthetic system in the inner iteration (Newton-MS) is constructed. Precisely, the nonlinear Boltzmann equation is solved by an outer modified Newton iteration, and each Newton correction equation is solved by an inner source iteration preconditioned by a macroscopic synthetic system. The linearized collision operator is used in the Newton correction equation to reduce the cost of collision evaluation, while the macroscopic correction accelerates the propagation of large-scale flow information when $\Kn$ is small.

For the modified Newton's method to solve the steady-state Boltzmann equation, it aims to seek the numerical solution to \eqref{eq:ss_bolt}
\begin{equation}
    \label{eq:ss_bolt}
    \mathcal{R}(f) := \bv \cdot \nabla_{\bx} f  - \frac{1}{\Kn}\mQ[f, f] = 0.
\end{equation}
A similar Newton method proposed in \cite{yin2025fast} is utilized here, where the distribution function is updated iteratively by 
\begin{equation} \label{eq:Newton_update}
    f^{(n+1)} = f^{(n)} - g^{(n)},
\end{equation}
where the correction $g^{(n)}$ is solved from the linearized equation
\begin{equation}
    \label{eq:Newton_eq}
    \bv \cdot \nabla_{\bx}  g^{(n)} - \frac{1}{\Kn}\Big(\mathcal{Q}[f^{(n)},g^{(n)}]+\mathcal{Q}[g^{(n)},f^{(n)}]\Big) = r^{(n)},
\end{equation}
with 
\begin{equation}
    \label{eq:r}
r^{(n)} = \mathcal{R}(f^{(n)})    
\end{equation}
the residual at the $n$-th Newton iterate. This Newton iteration \eqref{eq:Newton_update} to update the distribution function is called the outer Newton iteration. The computational cost can be extremely high when numerically solving \eqref{eq:Newton_eq} with an iterative method, regardless of which method is used, due to the high cost of computing the binary collision operator. To solve this, the quadratic collision term $\mathcal{Q}[f^{(n)},g^{(n)}]+\mathcal{Q}[g^{(n)},f^{(n)}]$ is approximated with a linearized collision operator around the local Maxwellian, with a similar method utilized in \cite{yin2025fast}. This strategy is conceptually similar to the quasi-Newton method by replacing the Jacobian matrix with a computationally efficient approximation. Therefore, the governing equation of the correction $g^{(n)}$ is reduced to
\begin{equation}
    \label{eq:Newton_eq2}
    \bv \cdot \nabla_{\bx}  g^{(n)} -\frac{1}{\Kn}\mathcal{L}^{(n)}[g^{(n)}] = r^{(n)},
\end{equation}
where the linearized operator $\mathcal{L}^{(n)}[\cdot]$ is defined around the local Maxwellian as 
\begin{equation}
    \label{eq:linear_coll}
    \mathcal{L}^{(n)}[\cdot] = \mathcal{Q}[\mathcal{M}^{(n)}, \cdot]+\mathcal{Q}[\cdot,\mathcal{M}^{(n)}].
\end{equation}
Here, $\mathcal{M}^{(n)}$ is the Maxwellian corresponding to $f^{(n)}$. It has been demonstrated in \cite{yin2025fast} that this modified Newton's method shows fast convergence rates and can greatly reduce the computational cost at the same time. 

Within each Newton's iteration, to solve the linear equation \eqref{eq:Newton_eq2}, the iterative method is also utilized, which is called the inner iteration. For example, the classical source iteration (SI) was employed in \cite{yin2025fast}, which updates the numerical solution by 
\begin{equation}
    \label{eq:SI}
    \bv \cdot \nabla_{\bx}  g^{(n,l+1)} + \frac{{\nu}}{\Kn} g^{(n,l+1)} = \frac{1}{\Kn}\mathcal{L}^{(n)}[g^{(n,l)}] + \frac{{\nu}}{\Kn}  g^{(n,l)} + r^{(n)},
\end{equation}
where $n$ stands for the index of the outer Newton iteration, and $l$ stands for the index of the inner source iteration. Here, the penalization terms involving $\nu$ are introduced to obtain a stable source iteration, where $\nu$ is usually chosen as the local collision frequency depending on the density and temperature. The resulting outer-inner solver is referred to as the modified Newton method combined with source iteration (Newton--SI). However, it is well known that SI leads to slow convergence when the Knudsen number is small \cite{su2020fast, wang2018comparative, wu2017fast, cai2025symmetric}.

Improving the computational efficiency of SI is full of challenges, and several numerical methods have been proposed to accelerate it, especially in the near-continuum regime. 
For example, a macroscopic synthetic preconditioner is constructed by solving a macroscopic system, and a similar macroscopic coupling principle is also proposed in the general synthetic iterative scheme (GSIS) \cite{su2020fast, su2020can}. 

The general idea of GSIS is to accelerate the convergence of macroscopic variables by solving an NSF-like system, where the stress tensor and the heat flux include higher-order contributions from the microscopic solution in the previous iteration. Precisely, the synthetic equations are the full nonlinear NSF system with high-order corrections in nonlinear GSIS \cite{zhu2021general}, while the linearized counterpart around a global equilibrium is adopted in \cite{su2020fast}. The main contribution of this work is to construct an efficient macroscopic synthetic preconditioner to solve \eqref{eq:SI}  by coupling the macroscopic governing equations.

\subsection{Macroscopic moment system for the Newton correction equation}
\label{sec:acc}
To improve the convergence of SI in the near-continuum regime, we follow the idea of adopting the macroscopic system to accelerate convergence, but the specific form of the macroscopic equations is different from GSIS. Compared with GSIS, the correction equation \eqref{eq:Newton_eq2} in this work is linearized around the local Maxwellian $\mathcal{M}^{(n)}(\bx, \bv)$, where the macroscopic variables, such as the density $\rho^{(n)}$ and macroscopic velocity $\bu^{(n)}$, vary in space and will only be refreshed at the outer Newton step when the distribution function is updated by \eqref{eq:Newton_update}. This may bring a different macroscopic synthetic preconditioner, which will be derived through the Chapman-Enskog expansion of the correction $g^{(n)}$ around $\mathcal{M}^{(n)}$. 

Before presenting the macroscopic synthetic preconditioner, the local equilibrium $\mathcal{E}^{(n)}$ of the correction $g^{(n)}$ is introduced. Defining the peculiar velocity $\boldsymbol{c}^{(n)} = \bv - \bu^{(n)}(\bx)$, the local Maxwellian $\mathcal{M}^{(n)}$ of the distribution function $f^{(n)}$ in \eqref{eq:linear_coll} can be rewritten as 
\begin{equation}
\label{eq:local_Max}
\mathcal{M}^{(n)}(\bx,\bv) = \frac{\rho^{(n)}(\bx)}{(2\pi T^{(n)}(\bx))^{3/2}} \exp \left( -\frac{|\boldsymbol{c}^{(n)}|^2}{2T^{(n)}(\bx)} \right),
\end{equation}
where $\rho^{(n)}(\bx)$, $\bu^{(n)}(\bx)$, and $T^{(n)}(\bx)$ denote the macroscopic moments of $f^{(n)}$ through \eqref{eq:mac_var}. During the inner iteration in one Newton step, these macroscopic variables are known and fixed. Then the kernel space of the linearized collision operator \eqref{eq:linear_coll} \cite{cai2015approximation} holds the form as 
\begin{equation}
\label{eq:ker_linear}
\operatorname{ker} \mathcal{L}^{(n)} = \{\varphi(\boldsymbol{c}^{(n)})^T \alpha(\bx) \mathcal{M}^{(n)}(\bx,\bv) \mid \alpha: \Omega \rightarrow \mathbb{R}^5\}, 
\end{equation}
where
\begin{equation}
    \label{eq:phi}
    \varphi(\boldsymbol{c}^{(n)}) = \left( 1, \boldsymbol{c}^{(n)}, \displaystyle \frac{|\boldsymbol{c}^{(n)}|^2 - 3T^{(n)}}{2} 
    \right)^T.
\end{equation}
Then, the local equilibrium $\mathcal{E}^{(n)}(\bx,\bv)$ of the correction $g^{(n)}$ is defined as the projection of $g^{(n)}$ onto $\operatorname{ker} \mathcal{L}^{(n)}$ 
\begin{equation}
    \label{eq:linear_Max}
    \mathcal{E}^{(n)}(\bx,\bv) \triangleq  \varphi(\boldsymbol{c}^{(n)})^T {m}^{(n)}(\bx) \mathcal{M}^{(n)}(\bx,\bv),
\end{equation}
where $\bm{m}^{(n)}$ is the vector of local macroscopic variables defined as 
\begin{equation}
    \label{eq:local_mac}
    {m}^{(n)}(\bx) = \begin{pmatrix}
        \displaystyle \frac{\delta \rho^{(n)}(\bx)}{\rho^{(n)}(\bx)}, &
        \displaystyle \frac{\delta \bu^{(n)}(\bx)}{T^{(n)}(\bx)}, &
        \displaystyle \frac{\delta T^{(n)}(\bx)}{(T^{(n)}(\bx))^2}
    \end{pmatrix}^T,
\end{equation}
with $\delta \rho$, $\delta \bu$, and $\delta T$ being the perturbations of macroscopic variables related to $\hatn{g}$
\begin{equation}
\label{eq:macro_g}
 {m}^{(n)}(\bx) = \int_{\bbR^3} \overline{\varphi}(\bm{c}^{(n)}) g^{(n)} \dd \bv, \qquad 
\overline{\varphi}(\bm{c}^{(n)}) = \left(\frac{1}{\hatn{\rho}}, \frac{\bm{c}^{(n)}}{\hatn{\rho}\hatn{T}}, \frac{|\hatn{c}|^2 - 3\hatn{T}}{3\hatn{\rho}(\hatn{T})^2}\right)^T.
\end{equation}
It is easy to verify that 
\begin{equation}
    \label{eq:conver_g}
   \int_{\mathbb{R}^3}  \phi(\bv) g^{(n)} \dd \bv =  \int_{\mathbb{R}^3}  \phi(\bv) \mathcal{E}^{(n)}  \dd \bv, \qquad \phi(\bv) = (1, \bv, |\bv|^2)^T. 
\end{equation}
Defining $\mathcal{P}$ as the orthogonal projection operator from any distribution function to $\operatorname{ker} \mathcal{L}$, it holds that 
\begin{equation}
    \label{eq:operator_P}
  \int_{\mathbb{R}^3}  \phi(\bv) (I - \mP) f \dd \bv = 0, \qquad \mP g^{(n)} = \mathcal{E}^{(n)}. 
\end{equation}
The detailed derivation of the local equilibrium of $g^{(n)}$ can be found in \cite{cai2015approximation}. 

In the following, the macroscopic system associated with the Newton correction will be derived with a Chapman-Enskog type expansion of the correction $g^{(n)}$ in terms of the Knudsen number, and the outer-iteration superscript ``$(n)$'' is omitted for simplicity. Thus, the linearized equation \eqref{eq:Newton_eq2} is reduced to 
\begin{equation}
    \label{eq:Kn_scaled}
    \bv \cdot \nabla_{\bx} g - \frac{1}{\Kn} \mathcal{L}[g] = r.
\end{equation}
Moreover, with the local Maxwelilan \eqref{eq:linear_Max}, the Chapman-Enskog ansatz of the correction $g$ has the form 
\begin{equation}
    \label{eq:ce_expan}
    g = \mathcal{E} + \Kn g_1 + \Kn^2 g_2,
\end{equation}
where $g_1$ represents the first-order Navier-Stokes correction, and $g_2$ denotes the high-order non-equilibrium remainder. With the definition of $\mP$, it holds that  
\begin{equation}
    \label{eq:g1_g2}
    \mathcal{P} g_1 = 0,  \qquad \mathcal{P} g_2 = 0.
\end{equation}
Substituting \eqref{eq:ce_expan} into \eqref{eq:Kn_scaled}, with $\mathcal{L}[\mathcal{E}] = 0$, we can obtain
\begin{equation} \label{eq:CE}
    \bv \cdot \nabla_{\bx} \mathcal{E} + \Kn\, \bv \cdot \nabla_{\bx} g_1 + \Kn^2 \bv \cdot \nabla_{\bx} g_2 - \mathcal{L}[g_1] - \Kn\, \mathcal{L}[g_2] = r.
\end{equation}
To obtain $g_1$, by comparing the terms in order $\mO(1)$, we obtain 
\begin{equation}
    \label{eq:order_1}
\mathcal{L}[g_1]  = \bv \cdot \nabla_{\bx} \mathcal{E}  - r.
\end{equation}
Modding out the kernel of $\mL$ by taking the operator $(I-\mP)$, then $g_1$ is derived as 
\begin{equation}
    \label{eq:ce_expan2}
    g_1 = \mathcal{L}^{-1}(I - \mathcal{P})\!\left(\bv \cdot \nabla_{\bx} \mathcal{E} - r\right),
\end{equation}
where $\mathcal{L}^{-1}$ is the pseudoinverse of $\mathcal{L}$. From \eqref{eq:operator_P}, it holds for $g_1$ and $g_2$ that 
\begin{equation}
    \label{eq:con_g1}
     \int_{\mathbb{R}^3}  \phi(\bv) g_1 \dd \bv =  \int_{\mathbb{R}^3}  \phi(\bv) g_2 \dd \bv = 0.
\end{equation}
Moreover, since $\mathcal{L}^{-1} (I-\mathcal{P}) = \mathcal{L}^{-1}$, the projection $(I-\mathcal{P})$ will be omitted below.

The macroscopic equations for $\delta \rho$, $\delta \bu$ and $\delta T$ are obtained by taking moments of \eqref{eq:CE} with respect to the collision invariants $\phi(\bv)$. With the conservation property of the linearized collision operator and the first-order correction \eqref{eq:ce_expan2}, the resulting macroscopic equations can be written as
\begin{subequations}
\label{eq:macro_expanded}
    \begin{align}
    \label{eq:macro_expanded_1}
       \mathcal{S}_1 &= \int_{\bbR^3} r \dd \bv, \\ 
               \label{eq:macro_expanded_2}
       \mathcal{S}_2 +  \Kn\, \nabla_{\bx} \cdot \boldsymbol{\sigma} &=  \int_{\bbR^3} \bv r \dd \bv - \Kn^2 \nabla_{\bx} \cdot \int_{\bbR^3} \bv\otimes \bv g_2 \dd \bv, \\
        \label{eq:macro_expanded_3}
        \mathcal{S}_3 + 2\Kn\, \nabla_{\bx} \cdot \left( \boldsymbol{\sigma} \bu + \bq\right) &= \int_{\bbR^3} |\bv|^2\, r \dd \bv - \Kn^2 \nabla_{\bx} \cdot \int_{\bbR^3} |\bv|^2 \bv\, g_2\dd \bv,
    \end{align}
\end{subequations}
with 
\begin{subequations}
\label{eq:eq_S}
    \begin{align}
       \mathcal{S}_1= &\nabla_{\bx} \cdot \left( \bu\, \delta\rho + \rho\, \delta\bu \right),  \\ 
      \mathcal{S}_2 = &\nabla_{\bx} \cdot \Big[ \delta\rho\, (\bu \otimes \bu + T \mathbf{I}) + \rho\, (\bu \otimes \delta \bu + \delta \bu \otimes \bu) + \rho\, \delta T\, \mathbf{I}\Big], \\
        \mathcal{S}_3=& \nabla_{\bx} \cdot \Big[ \delta\rho \bu (|\bu|^2 + 5T) + \rho\, \delta\bu\, (|\bu|^2 + 5T) + 2\rho \bu (\bu \cdot \delta\bu) + 5\rho\, \bu\, \delta T \Big].
    \end{align}
\end{subequations}
Here, the variables $\bm{\sigma}$ and $\bq$ are defined as 
\begin{equation}
\label{eq:sigma_q}
    \boldsymbol{\sigma} = \int_{\bbR^3} \left( \boldsymbol{c} \otimes \boldsymbol{c} - \frac{|\boldsymbol{c}|^2}{3}\mathbf{I} \right) g_1 \dd \bv, \qquad
    \boldsymbol{q}= \int_{\bbR^3} \frac{1}{2} |\boldsymbol{c}|^2 \boldsymbol{c}\, g_1 \dd \bv,
\end{equation}
with $\otimes$ denotes the tensor product of two vectors. The left-hand side of \eqref{eq:macro_expanded} has a similar form to the standard linearized Navier-Stokes-Fourier equations (NSF), which are approximated in the reference state $(\rho, \bu, T)$. But the variables $\bm{\sigma}$ and $\bq$ defined in \eqref{eq:sigma_q} depend not simply on $\delta \bu$ and $\delta T$ as the standard NSF, but also on the macroscopic variables $\rho, \bu$ and $\theta$ due to the non-homogeneous Maxwellian $\mM^{(n)}$ in the outer Newton iteration. To obtain $\bm{\sigma}$ and $\bm{q}$, they are split into two parts
\begin{equation}
\label{eq:split_s_q}
    \boldsymbol{\sigma} =  \boldsymbol{\sigma}^{(\nabla)} - \boldsymbol{\sigma}^{(r)}, \qquad
    \boldsymbol{q} =  \boldsymbol{q}^{(\nabla)} - \boldsymbol{q}^{(r)},
\end{equation}
where 
\begin{subequations} 
\begin{align}
\label{eq:sigma_q_g}
   \boldsymbol{\sigma}^{(\nabla)} &= \int_{\bbR^3} \left( \boldsymbol{c}\otimes \boldsymbol{c} - \frac{|\boldsymbol{c}|^2}{3} \mathbf{I} \right) \mathcal{L}^{-1} (\bv \cdot \nabla_{\bx} \mathcal{E}) \dd \bv, \qquad
    \boldsymbol{q}^{(\nabla)} = \int_{\bbR^3} \frac{|\boldsymbol{c}|^2}{2} \boldsymbol{c}\, \mathcal{L}^{-1}  (\bv \cdot \nabla_{\bx} \mathcal{E}) \dd \bv, \\
\label{eq:sigma_q_r}
      \boldsymbol{\sigma}^{(r)} &= \int_{\bbR^3} \left( \boldsymbol{c}\otimes \boldsymbol{c} - \frac{|\boldsymbol{c}|^2}{3} \mathbf{I} \right) \mathcal{L}^{-1} r \dd \bv, \qquad
    \boldsymbol{q}^{(r)} = \int_{\bbR^3} \frac{|\boldsymbol{c}|^2}{2} \boldsymbol{c}\, \mathcal{L}^{-1} r \dd \bv.
\end{align}  
\end{subequations}
Here, $ \boldsymbol{\sigma}^{(\nabla)}$ and $\boldsymbol{q}^{(\nabla)}$ are calculated as 
\begin{align}
\label{eq:sigma_q_full}
    \boldsymbol{\sigma}^{(\nabla)} &= \frac{C_{2,0}}{\nu}\,\mathbf{S}
        + \frac{C_{2,1}}{\nu}\,\mathbf{R}, \qquad 
    \boldsymbol{q}^{(\nabla)} = -\frac{1}{2\nu} \left({C_{1,1}}\,\,\boldsymbol{P}
        + {C_{1,2}}\,\,\boldsymbol{Q} \right),
\end{align}
where
{\small
\begin{subequations}
    \label{eq:eq_s}
\begin{align}
    \label{eq:S}
    \mathbf{S} &= \frac{1}{T}(\nabla_{\bx}\delta \bu)_{\mathrm{stf}} + \tf{\frac{\delta \bu}{T}\otimes\frac{\nabla_{\bx}\rho}{\rho}}
        + \left(\frac{\delta\rho}{\rho T} - \frac{5\delta T}{2T^2}
    \right)\tf{\nabla_{\bx}\bu} \notag \\
        &\quad - \frac{5}{2T^2}\tf{\delta\bu\otimes\nabla_{\bx}T}
        + \tf{\frac{\delta\bu}{T}\otimes(\bu\cdot\nabla_{\bx})\bu}, \\
    \label{eq:R}
    \mathbf{R} &= \frac{\rho\delta T}{3T}\,\tf{\nabla_{\bx}\bu}
        + \frac{\rho}{2T}\,\tf{\delta\bu\otimes\nabla_{\bx}T},
\end{align}
\end{subequations}
and
\begin{subequations}
\label{eq:eq_q}
\begin{align}
    \label{eq:P}
    \boldsymbol{P} &= \rho T\,\nabla_{\bx}\delta T
        + T\delta T\,\nabla_{\bx}\rho
        + T\delta\rho\,\nabla_{\bx} T
        - 5\rho\delta T\,\nabla_{\bx} T \notag \\       
        & \quad + \frac{2\rho T}{5}\bigl[(\delta\bu\cdot\nabla_{\bx})\bu
            + (\nabla_{\bx}\bu)^{\top}\delta\bu
            + (\nabla_{\bx}\cdot\bu)\,\delta\bu\bigr] \\ \notag
             &\quad + \rho\delta T\,(\bu\cdot\nabla_{\bx})\bu
        + \rho(\bu\cdot\nabla_{\bx} T)\,\delta\bu \\
    \label{eq:Q}
    \boldsymbol{Q} &= \frac{\rho\delta T}{2}\,\nabla_{\bx} T.
\end{align}
\end{subequations}
}
Here, $\tf{\mathbf{X}}$ denotes the symmetric and trace-free part of the matrix $\mathbf{X}$, defined by $(\mathbf{X} + \mathbf{X}^T)/2 - \operatorname{tr}(\mathbf{X})\mathbf{I}/3$. The detailed deduction is listed in App. \ref{sec:constitutive_derivation}. For the spatially homogeneous background where $\nabla_{\bx}\rho$,  $\nabla_{\bx}\bu$  and $\nabla_{\bx}T$ all equal zero, \eqref{eq:eq_s} and \eqref{eq:eq_q} reduce to the standard constitutive laws
\begin{equation}
\label{eq:sigma_q_frozen}
    \boldsymbol{\sigma}^{(\nabla)} = \frac{C_{2,0}}{\nu T} \tf{\nabla_{\bx} \delta\bu}, \qquad \bq^{(\nabla)} = -\frac{C_{1,1}}{2\nu}\rho T\, \nabla_{\bx} \delta T.
\end{equation}
Substituting \eqref{eq:split_s_q} in \eqref{eq:macro_expanded}, with \eqref{eq:sigma_q_full}, the macroscopic system \eqref{eq:macro_expanded} is reduced into
{\small
\begin{subequations}
\label{eq:macro_expanded_update}
    \begin{align}
         \label{eq:macro_expanded_update_1}
\mathcal{S}_1 &= \int_{\bbR^3} r \dd \bv, \\
         \label{eq:macro_expanded_update_2}
\mathcal{S}_2 + \Kn \nabla_{\bx}\cdot  \boldsymbol{\sigma}^{(\nabla)} &= \Kn \nabla_{\bx}\cdot  \boldsymbol{\sigma}^{(r)} +  \int_{\bbR^3} \bv r \dd \bv - \mathcal{G}_{2,1},\\
         \label{eq:macro_expanded_update_3}
\mathcal{S}_3 + 2\Kn\, \nabla_{\bx} \cdot \left( \boldsymbol{\sigma}^{(\nabla)} \bu + \bq^{(\nabla)}\right)
&= 2\Kn\, \nabla_{\bx} \cdot \left( \boldsymbol{\sigma}^{(r)} \bu + \bq^{(r)}\right) + \int_{\bbR^3} |\bv|^2\, r \dd \bv - \mathcal{G}_{2,2},
\end{align}
\end{subequations}
}
with
\begin{equation}
    \label{eq:G2}
    \mathcal{G}_{2,1} = \Kn^2 \nabla_{\bx} \cdot \int_{\bbR^3} \bv\otimes \bv g_2 \dd \bv,\qquad \mathcal{G}_{2,2} = \Kn^2 \nabla_{\bx} \cdot \int_{\bbR^3} |\bv|^2 \bv\, g_2\dd \bv
\end{equation}
where the left parts are the linear system with the unknown variables $\delta\rho$, $\delta \bu$ and $\delta T$, who depend on $\rho$, $\bu$ and $T$, which keep constant during the inner iteration. The right parts are unknown terms left to be closed, which will be introduced in detail in the following section. Then, the macroscopic system \eqref{eq:macro_expanded} will be completed, which form a linear system of $\delta \rho$, $\delta \bu$, and $\delta T$, while most of the coefficients are spatially dependent.

\paragraph{BGK-type approximation of the coefficients $C_{i,j}$}
For the linearized collision operator, to obtain the coefficients $C_{i,j}$ in \eqref{eq:sigma_q_full}, the operator $\mathcal{L}^{-1}$ should be calculated exactly, which is quite difficult \cite{cai2015approximation}. In the numerical implementation, the BGK model is utilized as a penalty term to obtain $C_{i,j}$. Precisely, \eqref{eq:Kn_scaled} is reformulated by inserting a simpler BGK collision operator:
\begin{equation} \label{eq:BGK_penalized}
  \bv \cdot \nabla_{\bx} g - \frac{1}{\Kn}\mL_{\rm BGK}[g] = r + \frac{1}{\Kn}\Big(\mathcal{L}[g] - \mL_{\rm BGK}[g]\Big), \qquad \mL_{\rm BGK} = \nu(\mathcal{E} - g).
\end{equation}
Thus, define 
\begin{equation}
    \label{eq:tilde_r}
    \tilde{r} = r + \Kn^{-1} \big(\mathcal{L}[g] - \mathcal{L}_{\rm BGK}[g]\big).
\end{equation}
It is easy to verify that utilizing $\tilde{r}$ instead of $r$ in \eqref{eq:BGK_penalized}, the deduction of the macroscopic synthetic system is the same as \eqref{eq:macro_expanded_update}, due to 
\begin{equation}
    \int \phi(\bv) r \dd\bv = \int \phi(\bv) \tilde{r} \dd\bv.
\end{equation}
Moreover, the pseudoinverse of the BGK operator has the form 
\begin{equation}
    \label{eq:BGK_inv}
    \mL_{\rm BGK}^{-1} =-\frac{1}{\nu} (I-\mathcal{P}),
\end{equation}
thus, the coefficients $C_{i,j}$ can be obtained explicitly as
\begin{equation}
    C_{2,0} = -\rho T^2,\qquad C_{2,1} = -7 T,\qquad C_{1,1} = 5,\qquad C_{1,2} = 70.
\end{equation}
The detailed deduction is presented in App. \ref{app:bgk_coefficients}.



\subsection{Closure and synthetic iteration}
\label{sec:newton_gsis_iter}
To obtain the right-side terms of \eqref{eq:macro_expanded_update}, we recall the numerical scheme of the inner iteration \eqref{eq:SI}. Omitting the superscript $(n)$ in the outer iteration, \eqref{eq:SI} is reduced into 
\begin{equation}
    \label{eq:SI_update}
     \bv \cdot \nabla_{\bx} g^{(l+1)} + \frac{{\nu}}{\Kn} g^{(l+1)} = \frac{1}{\Kn}\mathcal{L}^{(n)}[g^{(l)}] + \frac{{\nu}}{\Kn}  g^{(l)} + r,
\end{equation}
with $r$ defined in \eqref{eq:r}, which keeps constant in the inner iteration. Moreover, the background macroscopic variables $\rho,\bu$, and $T$ also remain constant. Let 
\begin{equation}
    \label{eq:g1_l}
      g_1^{(\nabla,l)}= \mathcal{L}^{-1}\left(\bv\cdot\nabla_{\bx}\mathcal{E}^{(l)}\right),
  \qquad
  g_1^{(r)} = g_1^{(r,l)} = \mathcal{L}^{-1}r, \qquad g_1^{(l)} = g_1^{(\nabla,l)} - g_1^{(r)},
\end{equation}
where $\mathcal{E}^{(l)}$ is the corresponding local equilibrium of $g^{(l)}$. Here, the superscript $(l)$ is omitted in $g_1^{(r, l)}$ due to the constancy of the residual $r$ \eqref{eq:r} in each inner iteration. With the Chapman-Enskog decomposition \eqref{eq:ce_expan} and \eqref{eq:ce_expan2}, it holds at the $l$-th inner iteration that
\begin{equation}
\label{eq:g2}
  \Kn^2 g_2^{(l)} =  g^{(l)} - \mathcal{E}^{(l)} - \Kn\,g_1^{(\nabla,l)} + \Kn g_1^{(r,l)}.
\end{equation}
 Moreover, due to \eqref{eq:con_g1}, we can derive that 
\begin{equation}
    \label{eq:G2_1}
  \int_{\bbR^3} \bv \otimes \bv g_2^{(l)} \dd\bv =  \int \left( \boldsymbol{c} \otimes \boldsymbol{c} - \frac{|\boldsymbol{c}|^2}{3}\mathbf{I} \right) g_2^{(l)} \dd\bv, \qquad  \int_{\bbR^3} |\bv|^2\bv  g_2^{(l)} \dd\bv =  \int_{\bbR^3} |\bm{c}|^2 \bm{c} g_2^{(l)} \dd\bv.
\end{equation}

Since the terms \eqref{eq:G2} are small variables, they are closed with variables obtained at the $l$-th inner iteration to close \eqref{eq:macro_expanded_update}. Precisely, substituting \eqref{eq:g2} into \eqref{eq:G2}, with \eqref{eq:G2_1}, \eqref{eq:sigma_q_r}, then \eqref{eq:macro_expanded_update} is reduced into 
{\small
\begin{subequations}
    \label{eq:macro_expanded_final}
\begin{align}
    \label{eq:macro_expanded_final_1}
    \mathcal{S}_1 &= \int_{\bbR^3} r \dd \bv, \\
         \label{eq:macro_expanded_final_2}
\mathcal{S}_2 + \Kn \nabla_{\bx}\cdot  \boldsymbol{\sigma}^{(\nabla)} &= \int_{\bbR^3} \bv r \dd \bv - \Kn \nabla_{\bx} \cdot \left(\boldsymbol{\sigma}^{(\mathrm{neq},l)} - \boldsymbol{\sigma}^{(\nabla,l)}\right),\\
         \label{eq:macro_expanded_final_3}
\mathcal{S}_3 + 2\Kn\, \nabla_{\bx} \cdot \left( \boldsymbol{\sigma}^{(\nabla)} \bu + \bq^{(\nabla)}\right)
&= \int_{\bbR^3} |\bv|^2\, r \dd \bv - 2\Kn \nabla_{\bx} \cdot \left[ \left( \boldsymbol{\sigma}^{(\mathrm{neq},l)} \bu + \bq^{(\mathrm{neq},l)} \right) - \left( \boldsymbol{\sigma}^{(\nabla,l)} \bu + \bq^{(\nabla,l)} \right) \right].
\end{align}
\end{subequations}
}
Here, the terms \eqref{eq:sigma_q_r} are canceled and there is no need to compute $\mathcal{L}^{-1}r$. $\boldsymbol{\sigma}^{(\mathrm{neq},l)}$ and $\boldsymbol q^{(\mathrm{neq},l)}$ are non-equilibrium moments defined as 
{\small
\begin{equation}
  \label{eq:sigma_q_neq}
\boldsymbol{\sigma}^{(\mathrm{neq},l)} =
  \frac{1}{\Kn}
  \int_{\mathbb R^3}
  \left(
  \bc\otimes\bc-\frac{|\bc|^2}{3}\mathbf I
  \right)
  \left( g^{(l)}- \mathcal{E}^{(l)}\right)\,\dd\bv,\qquad 
  \boldsymbol q^{(\mathrm{neq},l)}=
  \frac{1}{\Kn}
  \int_{\mathbb R^3}
  \frac{1}{2}|\bc|^2\bc\,\left( g^{(l)}- \mathcal{E}^{(l)}\right)\,\dd\bv.
\end{equation}
}
They are calculated directly by \eqref{eq:sigma_q_neq} in the $l$-th inner iteration since $g^{(l)}$ and $\mathcal{E}^{(l)}$ are already known for the moment. $\boldsymbol{\sigma}^{(\nabla, l)}$ and $\boldsymbol{q}^{(\nabla,l)}$ are defined in \eqref{eq:sigma_q_g} with $\mathcal{E}$ replaced by $\mathcal{E}^{(l)}$, and are calculated by \eqref{eq:sigma_q_full}. Besides, the terms $\int_{\bbR^3} r \dd \bv, \int_{\bbR^3} \bv r \dd \bv$ and $\int_{\bbR^3} |\bv|^2 r \dd \bv$ keep constant in the inner iteration, and only one calculation is needed. 

For now, the macroscopic moment system is completely derived. Next, we will introduce how it is adopted to accelerate inner iteration \eqref{eq:SI_update}. At each outer Newton iteration, the correction is initialized as $g^{(0)} = 0$. Then, at the $l$-th inner iteration, the following three steps are applied
\begin{itemize}
    \item[\bf Step 1:] Obtain the macroscopic variables $\boldsymbol m^{(l)}=(\delta\rho^{(l)}/\rho, \delta\bu^{(l)}/T,\delta T^{(l)}/T^2)$ related to $g^{(l)}$ by \eqref{eq:macro_g}, the non-equilibrium moments $\boldsymbol{\sigma}^{(\mathrm{neq},l)}, \boldsymbol q^{(\mathrm{neq},l)}$ by \eqref{eq:sigma_q_neq}, and  $\boldsymbol{\sigma}^{(\nabla,l)}, \boldsymbol{q}^{(\nabla,l)}$ by \eqref{eq:sigma_q_full}.
    \item[\bf Step 2:] Obtain the intermediate macroscopic variables 
    \begin{equation}
    \label{eq:mid_macro}
    {m}^{(l+1,\ast)} = \left(\frac{\delta \rho^{(l+1, \ast)}}{ \rho}, \frac{\delta \bu^{(l+1,\ast)}}  {T}, \frac{\delta T^{(l+1,\ast)}}{ T^2}\right)^T 
    \end{equation}
    by solving the closed macroscopic system \eqref{eq:macro_expanded_final}, once  $\boldsymbol{\sigma}^{(\mathrm{neq},l)}$, $\boldsymbol{\sigma}^{(\nabla,l)}$ and $\boldsymbol{q}^{(\mathrm{neq},l)}$, $\boldsymbol{q}^{(\nabla,l)}$ are known. Then, reconstruct the local equilibrium part $\mathcal{E}^{(l+1,\ast)}$ from ${m}^{(l+1,\ast)}$.
     \item[\bf Step 3:] Obtain the Newton correction $g^{(l+1)}$ by the modified source iteration
     \begin{equation} 
     \label{eq:acc_SI}
  \bv\cdot\nabla_{\bx}g^{(l+1)} + \frac{{\nu}}{\Kn} \,g^{(l+1)}
  = \frac{1}{\Kn}\mathcal{L}\!\left[g^{(l)}\right] + \frac{{\nu}}{\Kn}\,g^{(l)}
  + \alpha^{(l)} \frac{{\nu}}{\Kn}\!\left(\mathcal{E}^{(l+1,\ast)} - \mathcal{E}^{(l)}\right) + r,
\end{equation}
where $\alpha^{(l)} \in [0,1]$ is a relaxation parameter.  
\begin{remark}
    With the macroscopic correction 
\begin{equation}
    \label{eq:macro_corr}
   \alpha \frac{{\nu}}{\Kn}(\mathcal{E}^{(l+1,\ast)} - \mathcal{E}^{(l)}),
\end{equation}
the information of the updated macroscopic variables is injected into the microscopic inner iteration, and therefore, the convergence of the macroscopic component in the correction $g^{(l)}$ is accelerated, especially when $\Kn$ is small. A similar acceleration method can also be found in \cite{su2020can}.
\end{remark}
\end{itemize}

To distinguish from Newton-SI, we call this Newton method accelerated by the macroscopic synthetic system in the inner iteration as Newton-MS. Though the convergence of the inner iteration \eqref{eq:acc_SI} can be accelerated with the macroscopic synthetic system \eqref{eq:macro_expanded_final}, directly solving \eqref{eq:macro_expanded_final} is still expensive. In the following section, the discrete form of \eqref{eq:macro_expanded_final} is directly deduced in the framework of the discontinuous Galerkin (DG) method, which is reduced into a linear system of $g^{(n)}$ and can be directly solved.

\section{Discrete form of Newton-MS}
\label{sec:num_dis}
In this section, the fully discrete form of Newton-MS will be introduced, and the whole process can be treated as deducing the discrete form of the macroscopic synthetic system \eqref{eq:macro_expanded_final} in the framework of DG. We begin from the discretization of the spatial space where the DG method is adopted. The physical domain $\Omega \subset \mathbb{R}^d$ ($d=1,2$) is partitioned into $N_{\mathrm{el}}$ elements, and on each element, the numerical solution is approximated by polynomials of degree $N_p$. For the microscopic velocity space, as in Sec. \ref{sec:fsm}, the microscopic velocity space is first truncated to a finite domain as $[-L, L]^3$, and then the uniform mesh with number $2N$ is utilized to discretize this microscopic velocity space. Thus, the discrete velocity set is $\{\boldsymbol{v}_k\}_{k=1}^{(2N)^3}$, where $\boldsymbol{v}_k=(v_{1k},v_{2k},v_{3k})^T$, and the corresponding quadrature weight is $\Delta v=(L/N)^3$. For simplicity, we introduce the total degrees of freedom $N_{t}$, those in the spatial space $N_{s}$ and those in the microscopic space $N_{v}$ as 
\begin{equation}
    \label{eq:dof}
    N_t = N_s \times N_v, \qquad N_{s} = N_{\mathrm{el}}(N_p+1)^d, \qquad N_{v} = (2N)^3.
\end{equation}
With this discretization, let $\boldsymbol{f}^{(n)}, \boldsymbol{g}^{(n)} \in \mathbb{R}^{N_t}$ be the global numerical solution and the global discrete correction, and their entries are ordered first by the spatial index. In particular,
\begin{equation}
   \boldsymbol{f}^{(n)} =\left((\boldsymbol{f}_1^{(n)})^T, \cdots, (\boldsymbol{f}_{N_{s}}^{(n)})^T \right)^T,\quad  \boldsymbol{g}^{(n)} =\left((\boldsymbol{g}_1^{(n)})^T, \cdots, (\boldsymbol{g}_{N_{s}}^{(n)})^T \right)^T,\qquad \boldsymbol{f}_i^{(n)}, \boldsymbol{g}_i^{(n)}\in\mathbb{R}^{N_v},
\end{equation}
where $\boldsymbol{f}_i^{(n)}, \boldsymbol{g}_i^{(n)}$ represent the distribution function $\bm{f}$, and the correction $\boldsymbol{g}$ at the $i$-th mesh, respectively. Thus, the update in each Newton step \eqref{eq:Newton_update} becomes 
\begin{equation}
    \label{eq:dis_Newton_update}
    \boldsymbol{f}^{(n+1)} = \boldsymbol{f}^{(n)} - \boldsymbol{g}^{(n)}. 
\end{equation}
Let $\mathbf{T}$ denote the DG discretization of the transport operator $\bv\cdot\nabla_{\bx}$. With the assumption of the upwind numerical flux with a homogeneous inflow condition, the operator $\mathbf T$ is reduced to a matrix as $\mathbf{T}\in\mathbb{R}^{N_{t}\times N_{t}}$. Moreover, it is block-tridiagonal in 1D and block-pentadiagonal in 2D on a structured mesh, with each block of size $(N_v (N_p+1)^d) \times (N_v (N_p+1)^d)$, due to the coupling with adjacent spatial elements through the upwind flux. In this case, the physical boundary condition is incorporated explicitly through a boundary source term. Thus, the discrete form \eqref{eq:Newton_eq2} has the form as 
\begin{equation} \label{eq:fully_discrete}
    \mathbf{T}\boldsymbol{g}^{(n)} - \frac{1}{\Kn}\mathbf{L}^{(n)}\boldsymbol{g}^{(n)} = \boldsymbol{r}^{(n)}+\boldsymbol{b}_{\rm bc}^{(n)},
\end{equation} 
where $\mathbf{L}^{(n)} \in \mathbb{R}^{N_t \times N_t}$ is the discrete form of the linearized collision operator \eqref{eq:linear_coll}, $\boldsymbol{r}^{(n)}$ is the discrete outer residual, and $\boldsymbol{b}_{\rm bc}^{(n)}$ is the boundary source term, depending on $\boldsymbol{g}^{(n)}$. Since the linear collision operator $\mathcal{L}^{(n)}$ is a local operator, $\mathbf{L}^{(n)}$ is block-diagonal
\begin{equation}
\label{eq:diag_L}
    \mathbf{L}^{(n)} = \operatorname{diag}\left\{
    \mathbf{L}_1^{(n)},\, \mathbf{L}_2^{(n)},\, \ldots,\, \mathbf{L}_{N_{s}}^{(n)}
    \right\},
\end{equation}
with each block $\mathbf{L}_i^{(n)}\in\mathbb{R}^{N_v\times N_v}$. In the following, the outer iteration superscript ``$(n)$'' will be omitted for brevity. 

To deduce the discrete form of the macroscopic synthetic system \eqref{eq:macro_expanded_final}, the discrete form of the deduction process as in Sec. \ref{sec:method} is proposed. Precisely, with the Chapman-Enskog ansatz \eqref{eq:ce_expan}, the discrete correction $\boldsymbol{g}$ is expanded as
\begin{equation}
    \label{eq:dis_ce}
    \boldsymbol{g} = \boldsymbol{\mathcal{E}} + \Kn \tilde{\boldsymbol{g}}_1 + \Kn^2\tilde{\boldsymbol{g}}_2, \qquad \tilde{\bm{g}}_1, \tilde{\bm{g}}_2 \in \bbR^{N_t},
\end{equation}
where $\boldsymbol{\mathcal{E}}$ is the discrete local equilibrium of $\boldsymbol{g}$, $\tilde{\boldsymbol{g}}_1$ is the discrete first-order Navier–Stokes correction, and $\tilde{\boldsymbol{g}}_2$ represents the discrete high-order non-equilibrium remainder.
Substituting \eqref{eq:dis_ce} into the discrete Newton correction equation \eqref{eq:fully_discrete}, and using $\mathbf{L}  \boldsymbol{\mathcal{E}}=0$, we obtain
\begin{equation} \label{eq:dis_CE_eq}
    \mathbf{T} \boldsymbol{\mathcal{E}} + \Kn\, \mathbf{T} \tbg_1 + \Kn^2 \mathbf{T} \tbg_2 - \mathbf{L}\tilde{\boldsymbol{g}}_1 - \Kn\, \mathbf{L}\tilde{\boldsymbol{g}}_2 = \boldsymbol{r}+\boldsymbol{b}_{\rm bc}.
\end{equation}
To obtain the discrete macroscopic equation related to \eqref{eq:macro_expanded}, the discrete collision invariant operator $ \Phi \in \mathbb{R}^{5N_{s}\times N_{t}}$ is introduced as 
\begin{equation}
\label{eq:dis_V}
    \Phi = \operatorname{diag}\{\bar{\Phi}, 
    \ldots, \bar{\Phi}\}, \qquad  \bar{\Phi} = \Delta v \left(
    \phi(\bv_1), \cdots, \phi(\bv_{N_v})
    \right)
    \in \mathbb{R}^{5\times N_v},
\end{equation}
where $\phi(\bv)$ is the collision invariants \eqref{eq:conver_g}. Multiplying \eqref{eq:dis_CE_eq} by $\Phi$ from the left, and with 
the discrete conservation property
\begin{equation}
    \label{eq:dis_conver}
\Phi\mathbf{L}=0,
\end{equation}
we can derive the discrete form of \eqref{eq:macro_expanded} as 
\begin{equation}
    \label{eq:dis_mac_eq1}
    \Phi \mathbf T\left( \boldsymbol{\mathcal{E}}+ \Kn \tbg_1 + \Kn^2\tbg_2\right) = \Phi\bm{r} + \Phi\boldsymbol b_{\rm bc}.
\end{equation}
Then, a similar closure is displayed to obtain the closed discrete macroscopic synthetic system. We first introduce two block-diagonal operators $\mathbf{S}$ and $\Gamma$ as 
\begin{equation}
\label{eq:S_G}
    \mathbf{S} = \operatorname{diag}\{\mathbf{S}_1,\ldots,\mathbf{S}_{N_s}\},
    \qquad
    \Gamma = \mathrm{diag}\{\Gamma_1,\ldots,\Gamma_{N_s}\},
\end{equation}
with $\bfS_i$ and $\Gamma_i$ as 
\begin{align}
\label{eq:SG}
    \mathbf{S}_i &= \Delta v\left(s_i(\bv_1),\cdots, s_i(\bv_{N_v})\right)\in \mathbb{R}^{5\times N_v}, \qquad 
    s_i(\bv_k) = \overline{\varphi}(\bm{c}_{ki}),
   \\
    \Gamma_i &= \left(\gamma_i(\bv_1), \cdots, \gamma_i(\bv_{N_v})\right)^T \in \mathbb{R}^{N_v\times 5}, \qquad\hspace{8pt} \gamma_i(\bv_k) =  \mathcal M_i(\bv_k)\varphi(\bc_{ki}). 
\end{align}
Here, $\bc_{ki} = \bv_k - \bu_i$ is the peculiar velocity, $\overline{\varphi}$ and $\varphi$ are defined in \eqref{eq:macro_g}, and \eqref{eq:phi}, respectively. $\mathcal{M}_i$ is the Maxwellian \eqref{eq:Maxwell}, whose corresponding macroscopic variables are $(\rho_i, \bu_i, T_i)$ at the $i$-th mesh. It is easy to derive that 
\begin{equation}
    \label{eq:oper_S}
    m_i = \bfS_i \bm{g}_i, \qquad \bm{\mathcal{E}}_i =  \Gamma_i m_i, \qquad i = 1, \cdots, N_s,
 \end{equation}
 It means that the operator $\bfS$ maps the discrete correction $\bm{g}$ to the local macroscopic variables $m_i$ \eqref{eq:local_mac}, and $\Gamma$ reconstructs local equilibrium from the local macroscopic variables. Moreover, we can deduce that $\Gamma\bfS$ is equal to the discrete form of the orthogonal projection operator $\mathcal{P}$ \eqref{eq:operator_P}. With the definition \eqref{eq:g1_l}, it holds that 
 \begin{equation}
     \label{eq:dis_g1}
      \tbg_1^{(l)} = \tbg_1^{(\nabla, l)} - \tbg_1^{(r)}, \qquad \tbg_1^{(\nabla, l)}=\mathbf{L}^{-1} \mathbf T\Gamma\bm{m}^{(l)},
 \end{equation}
where $\bm{m} = (m_1^T, \cdots, m_{N_s}^T)^T\in \bbR^{5N_s}$, and $\mathbf{L}^{-1}$ is the discrete form of the operator $\mathcal{L}^{-1}$. In this case, with the discrete Chapman-Enskog expansion \eqref{eq:g2}, \eqref{eq:dis_ce} and \eqref{eq:oper_S}, the high order $\tbg_2^{(l)}$ is obtained as 
\begin{equation}
\label{eq:dis_g2}
    \Kn^2 \tbg_2^{(l)} = \boldsymbol{g}^{(l)} - \Gamma \bm{m}^{(l)} - \Kn\tbg_1^{(\nabla,l)}+ \Kn \tbg_1^{(r)} \triangleq \hat{\bm{g}}^{(l)} +\Kn \tbg_1^{(r)}.
\end{equation}
Substituting \eqref{eq:oper_S}, \eqref{eq:dis_g1} and \eqref{eq:dis_g2} into the discrete moment equation \eqref{eq:dis_mac_eq1}, the residual-driven term $\boldsymbol{g}_1^{(r)}$ is canceled as in the continuous case \eqref{eq:macro_expanded_final}. Then, the discrete macroscopic system for the local macroscopic variables $\bm{m}^{(l+1,\ast)}$ in the inner iteration is driven as 
\begin{equation}
\label{eq:dis_mac_eq2}
    \Phi \mathbf T\left(\Gamma \bm{m}^{(l+1,\ast)} + \Kn\,\tbg_1^{(\nabla, l+1,\ast)} + \hat{\boldsymbol g}^{(l)}\right)=\Phi\bm{r}+ \Phi \boldsymbol b_{\rm bc}^{(l)}.
\end{equation}
The macroscopic system \eqref{eq:dis_mac_eq2} is the discrete form of the continuous closed system \eqref{eq:macro_expanded_final} and can be rewritten as
\begin{equation}
\label{eq:dis_mac_eq3}
    \left(\Psi_{\mathrm{E}}+\Kn\,\Psi_{\mathrm{NS}}\right)\bm{m}^{(l+1,\ast)}=\boldsymbol b^{(l)},
\end{equation}
where
\begin{equation}
\label{eq:matrix_K}
    \Psi_{\mathrm{E}}=\Phi\mathbf T\Gamma,\qquad\Psi_{\mathrm{NS}}=\Phi\mathbf T \mathbf L^{-1} \mathbf T\Gamma, \qquad 
    \boldsymbol b^{(l)}=-\Phi\mathbf T\hat{\boldsymbol g}^{(l)}+\Phi\bm{r}+\Phi\boldsymbol b_{\rm bc}^{(l)}.
\end{equation}
Here, $\Psi_{\mathrm{E}}$ represents the contribution of the equilibrium as \eqref{eq:eq_S}, while  $\Psi_{\mathrm{NS}}$ corresponds to the discrete first-order Chapman-Enskog correction, or the terms $(\cdot)^{(\nabla)}$ on the left side of \eqref{eq:macro_expanded_final}. 
\begin{remark}
    Due to the complex form of the discrete pseudoinverse $\mathbf L^{-1}$, the BGK-type pseudoinverse \eqref{eq:BGK_inv} is adopted here to approximate that of the linearized collision operator. Thus, $\Psi_{\rm NS}$ is reduced into 
\begin{equation}
    \label{eq:final_Psi}
   \Psi_{\mathrm{NS}} \approx -\frac{1}{\bm{\nu}}\Phi\mathbf T (\mathbf I - \Gamma \bfS) \mathbf T\Gamma,
\end{equation}
where $\bm{\nu} \in \bbR^{N_t}$ is the local collision frequency. 
\end{remark}
Note that $\Psi_{\rm E}$ and $\Psi_{\rm NS}$ depend only on the outer iteration variables and remain constant in the inner iterations and their reconstruction is based on applications of the sparse transport operator $\mathbf T$ on the block-diagonal matrix $\Gamma$, and the applications of the block-diagonal matrix $\Phi$, $\bfS$ and $\Gamma$.

 Finally, the discrete form of the augmented source iteration \eqref{eq:acc_SI} is 
\begin{equation}
\label{eq:dis_inner_si}
    \left(\mathbf T + \frac{\boldsymbol{\nu}}{\Kn}\right)\boldsymbol g^{(l+1)} = \frac{1}{\Kn}\mathbf L\boldsymbol{g}^{(l)} + \frac{\boldsymbol{\nu}}{\Kn}\boldsymbol g^{(l)} + \alpha \frac{\boldsymbol{\nu}}{\Kn}\Gamma\bigl(\bm{m}^{(l+1,\ast)}-\bm{m}^{(l)}\bigr) + \bm{r} + \boldsymbol b_{\rm bc}^{(l)}.
\end{equation}

For the computational complexity, the one-time setup within each Newton step consists of assembling the synthetic macroscopic operators $\Psi_{\mathrm{E}}$ and $\Psi_{\mathrm{NS}}$, whose construction is approximated as $\mathcal{O}(N_s N_p^dN^3)$.
For one inner iteration, the dominant cost comes from the implement of the linearized collision operator $\mathbf L\boldsymbol g^{(l)}$, whose complexity is $\mathcal{O}(N_sN^4\log N)$ with the FFT-based algorithm in Sec.~\ref{sec:fsm}. The computational cost brought by the transport operator $\mathbf T$ is approximated as $\mathcal{O}(N_sN_p^{d}N^3)$, while the numerical cost of solving the macroscopic synthetic system is $\mathcal{O}(N_s)$. Therefore, compared with Newton-SI, Newton-MS has essentially the same leading-order cost in the microscopic iteration. For the completeness of Newton-MS, the algorithm is listed in Alg. \ref{alg:discrete_newton_gsis}.

\begin{algorithm}[htbp]
\caption{Discrete Newton-MS}
\label{alg:discrete_newton_gsis}
\begin{algorithmic}[1]
\State \textbf{Input:} Initial distribution $\boldsymbol{f}^{(0)}$, 
Knudsen number $\Kn$, relaxation parameter $\alpha$
\State \textbf{Output:} Steady-state distribution $\boldsymbol{f}$
\State $n\gets 0$
\Repeat \Comment{Outer Newton iteration}
    \State Compute macroscopic moments $(\rho^{(n)},\bu^{(n)},T^{(n)})$ 
    from $\boldsymbol{f}^{(n)}$ through $\Phi$;
    \State Construct the local Maxwellian $\mathcal{M}^{(n)}$ and compute 
    the discrete residual $\boldsymbol{r}^{(n)}$;
    
    \State Build $\mathbf{S}$, $\Gamma$ from 
    $(\rho^{(n)},\bu^{(n)},T^{(n)})$ by \eqref{eq:SG}, 
    and assemble $\Psi_{\mathrm{E}}, \Psi_{\mathrm{NS}}$ by \eqref{eq:matrix_K};
    \State Initialize $\boldsymbol{g}^{(0)}\gets\boldsymbol{0}$, 
    $\bm{m}^{(0)}\gets\boldsymbol{0}$, $l\gets 0$;
    \Repeat \Comment{Inner accelerated source iteration}
        \State Compute 
        $\boldsymbol{\mathcal{E}}^{(l)} = \Gamma\bm{m}^{(l)}$
        and $\boldsymbol g_1^{(\nabla,l)}=-\frac{1}{\boldsymbol{\nu}}(\mathbf I-\Gamma\mathbf S) \mathbf T\Gamma\bm{m}^{(l)}$;
        \State Compute the remainder $\hat{\boldsymbol g}^{(l)}$ by \eqref{eq:dis_g2};
        \State Assemble the right-hand side  $\boldsymbol b^{(l)}$ by \eqref{eq:matrix_K};
        \State Solve the macroscopic system \eqref{eq:dis_mac_eq3} for $\bm{m}^{(l+1,\ast)}$;
        \State Solve the augmented source iteration for $\boldsymbol g^{(l+1)}$ by \eqref{eq:dis_inner_si};
        \State Set $\bm{m}^{(l+1)} = \mathbf{S}\boldsymbol{g}^{(l+1)}$;
        \State $l\gets l+1$;
    \Until{inner convergence}
    \State Set $\boldsymbol{g}^{(n)} = \boldsymbol{g}^{(l)}$ and update the distribution:
    $\boldsymbol{f}^{(n+1)} = \boldsymbol{f}^{(n)} - \boldsymbol{g}^{(n)}$;
    \State $n\gets n+1$;
\Until{outer convergence}
\State \Return $\boldsymbol{f}^{(n)}$
\end{algorithmic}
\end{algorithm}
\section{Numerical results}
\label{sec:num}

In this section, several numerical experiments are presented to validate the accuracy and efficiency of the proposed Newton-MS scheme. The main comparison is made with the standard Newton source iteration (Newton-SI), focusing on convergence behavior and computational cost over a range of Knudsen numbers, especially for the small Knudsen number. The collision operator is evaluated by the fast Fourier spectral method described in Sec. \ref{sec:fsm}, and the transport term is discretized by the discontinuous Galerkin (DG) method described in Sec. \ref{sec:num_dis}. 

The stop criteria for the outer Newton iteration  \eqref{eq:Newton_update} is that the residual of the nonlinear steady Boltzmann equation satisfies
\begin{equation}
\label{eq:residual}
    R_{\mathrm{out}}(n) = \left(\int_{\Omega}\int_{[-L,L]^3}\left|r^{(n)} \right|^2 
    \dd \bv \dd \bx\right)^{1/2} < \epsilon_{\rm out}.
\end{equation}
Moreover, the stop criteria for the inner iteration \eqref{eq:SI_update} are 
\begin{equation}
\label{eq:resi_inner}
    R_{\mathrm{in}}(n,l)=\left(\int_{\Omega}\int_{[-L,L]^3}\left|
    \mathcal R_{\mathrm{in}}^{(n,l)}\right|^2\dd \bv \dd \bx\right)^{1/2}< \epsilon_{\rm in,1}, 
        \quad\text{or}\quad
\frac{R_{\mathrm{in}}(n,l)}{R_{\mathrm{out}}(n)}<\epsilon_{\rm in, 2},
\end{equation}
where the residual of the inner iteration $\mathcal R_{\mathrm{in}}^{(n,l)}$ is defined as 
\begin{equation}
    \label{eq:inner_residual_def}
    \mathcal R_{\mathrm{in}}^{(n,l)} = \bv \cdot \nabla_{\bx} g^{(n,l)} - \frac{1}{\Kn}\mathcal{L}^{(n)}[g^{(n,l)}] - r^{(n)}.
\end{equation}
Here, the relative error in the stop criterion \eqref{eq:resi_inner} prevents unnecessary over-iteration of the linearized correction equation for $g$ when a highly accurate inner solution is not required for the outer Newton iteration. In the simulation, the related parameters $\epsilon_{\rm out}$, $\epsilon_{\rm in, 1}$ and $\epsilon_{\rm in, 2}$ are problem dependent, and unless otherwise specified, they are set as \eqref{eq:ep} in this work. 
\begin{equation}
    \label{eq:ep}
    \epsilon_{\rm out} = 10^{-5}, \qquad \epsilon_{\rm in, 1} = 10^{-6}, \qquad \epsilon_{\rm in, 2} = 10^{-2}.
\end{equation}

Moreover, the macroscopic correction 
\begin{equation}
    \label{eq:macro_cor}
    \alpha^{(l)} \frac{\nu}{\Kn}\left(\mathcal{E}^{(l+1, \ast)} - \mathcal{E}^{(l)}\right) 
\end{equation}
in \eqref{eq:acc_SI} mainly speeds up the convergence of macroscopic variables, including density, macroscopic velocity,  and temperature. Therefore, the inner residual drops quickly in the first few iterations, but it may stop decreasing once these macroscopic variables have converged, and the residual mainly comes from the non-equilibrium parts. In this case, the macroscopic correction \eqref{eq:macro_cor} is no longer needed to achieve convergence of the correction $g$. Thus, the relaxation parameter $\alpha^{(l)}$ in \eqref{eq:macro_cor} is chosen as
\begin{equation}
    \alpha^{(n,l)} =
        \begin{cases}
            \alpha_0, & l < l_{\mathrm{sw}}^{(n)},\\
            0, & l \geqslant l_{\mathrm{sw}}^{(n)},
        \end{cases}
    \label{eq:relax_para}
\end{equation}
where $l_{\mathrm{sw}}^{(n)}$ is the switching index at the $n$-th
Newton step, which is determined by the following residual-decay criterion
\begin{equation}
    l_{\mathrm{sw}}^{(n)}
    =
    \min\left\{
    l\geqslant p+1:
    \frac{R_{\mathrm{in}}(n,l)}
    {R_{\mathrm{in}}(n,l-j)}
    >
    \eta_{\mathrm{sw}},
    \quad j=1,\ldots,p
    \right\}.
    \label{eq:switch_criterion}
\end{equation}
Here, $p$ chosen as $p = 3$ is the monitoring window length and $\eta_{\mathrm{sw}}=0.9$ is a threshold to measure the rate of convergence. In the following tests, $\alpha_0=0.4$ and $\alpha_0 = 1$ are utilized for the one-dimensional Fourier and Couette flow problems, and for the two-dimensional cavity-flow problems, respectively.

For the truncation region of the microscopic velocity space, $L$ is set as 
\begin{equation}
\label{eq:L}
    L=\frac{3+\sqrt{2}}{2}R,\qquad R=4.
\end{equation}
Besides, when the fast Fourier spectral method is utilized to solve the linearized collision term, the total mass conservation can not be conserved \cite{yin2025fast}. Therefore, a post-processing is added to keep mass conservation. Precisely, after each outer Newton iteration, the distribution function is rescaled by 
\begin{equation}
  \label{eq:mass_correction}
  f^{(n+1)} \leftarrow f^{(n+1)}\cdot
  \frac{m_{\text{tot}}^{(0)}}{m_{\text{tot}}^{(n+1)}},
  \qquad
  m_{\text{tot}}^{(n+1)} = \int_\Omega\int_{[-L,L]^3} f^{(n+1)}\,\mathrm{d}\bv\,\mathrm{d}\bx,
\end{equation}
where $m_{\text{tot}}^{(0)}$ is the initial total mass.

In the simulation, the efficiency comparison between Newton-MS, Newton-SI and GSIS is displayed. For all three methods, they are all implemented with the same spatial and velocity discretization. For Newton-MS and Newton-SI, the outer iteration is the same nonlinear Newton iteration \eqref{eq:Newton_update}, but the macroscopic moment equations are utilized to accelerate the convergence \eqref{eq:acc_SI} of the inner iteration in Newton-MS, while \eqref{eq:SI} is directly adopted in Newton-SI. For GSIS, only the outer iteration with macroscopic synthetic acceleration is adopted, and no inner iteration is needed. 

For all three methods, to obtain the nonlinear outer residual $r^{(n)}$ \eqref{eq:r}, although the quadratic collision term only needs to be calculated once in each outer iteration, this is still quite expensive, and the CPU time to obtain $r^{(n)}$ for all $(n)$ is labeled $T_{\rm out}$. For Newton-MS and Newton-SI, to obtain the linear inner residual $\mathcal{R}_{\rm in}^{(n,l)}$ \eqref{eq:inner_residual_def},  only the linear collision term is calculated, but it needs to be calculated once for each inner iteration, which will also be quite expensive if the number of inner iterations is large. Here, the CPU time to obtain $\mathcal{R}_{\rm in}^{(n,l)}$ for all $(n,l)$ is labeled $T_{\rm in}$. Moreover, the CPU time for solving the macroscopic synthetic equations for all $(l,n)$ is labeled $T_{\rm m}$ with $T_{\rm tol}$ the total CPU time. 

\subsection{1D Fourier flow problem}
\label{sec:fourier}
We first consider the classical one-dimensional Fourier flow problem. The scenario consists of two parallel, stationary plates at $x=0$ and $x=1$, with wall temperatures $T_L=1$ and $T_R=1.2$, respectively. The discretization of the spatial domain is $N_{\rm el} = 40$ with the degree of polynomial $N_p = 2$. For the discretization of the microscopic velocity space, the number of the Fourier modes is set as $N = 24$, corresponding to the degree of freedom $2N=48$ in each direction of the microscopic velocity space. The fully diffuse reflection boundary condition \cite{yin2025fast} is imposed on both walls. The similar Fourier flow problem is also studied in \cite{jaiswal2019discontinuous}.

\begin{figure}[!hptb]
\centering
\subfigure[$\rho$]{
\includegraphics[width = 0.45\textwidth]{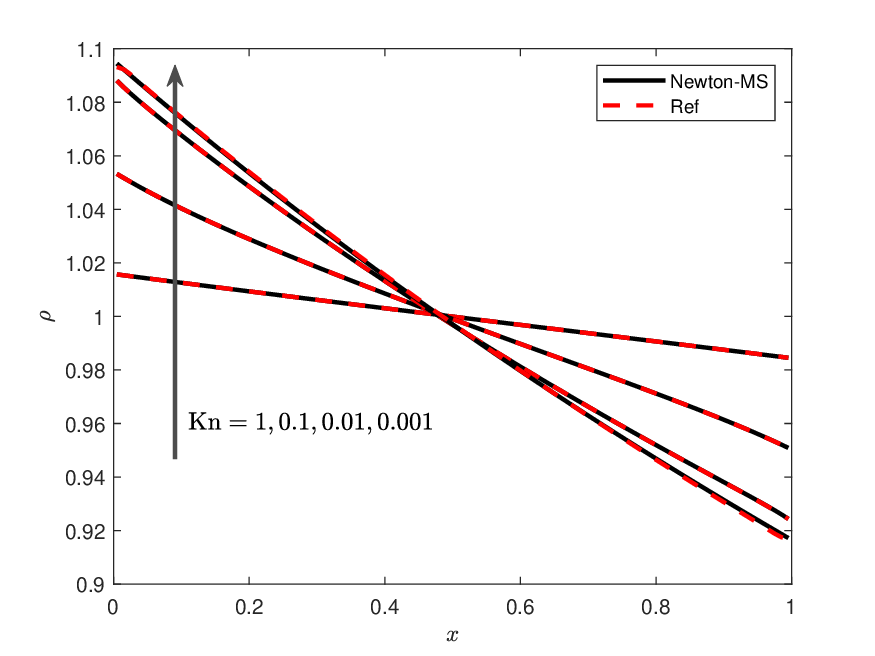}
}\qquad 
\subfigure[$T$]{
\includegraphics[width = 0.45\textwidth]{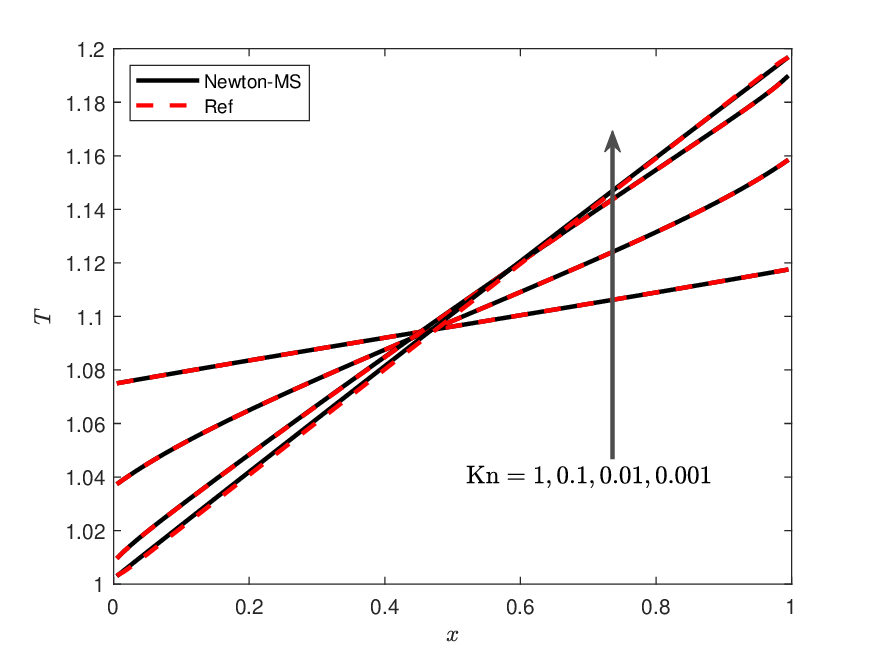}
}
\caption{(1D Fourier flow problem in Sec. \ref{sec:fourier}) Steady state numerical solution of the density $\rho$ and temperature $T$ for $\Kn=1,0.1,0.01$, and $0.001$. Here, the solid black lines are the numerical solution by Newton-MS, and the red dashed lines are those of the reference solution. }
\label{fig:fourier}
\end{figure}

The numerical solution of the density $\rho$, temperature $T$ at the steady state for $\Kn = 1, 0.1, 0.01$ and $0.001$ by Newton-MS is shown in Fig. \ref{fig:fourier}, where the reference solution for $\Kn = 0.001$ is obtained by directly solving the Navier-Stokes equation with temperature slip boundary condition imposed and those for other $\Kn$ are obtained by Newton-SI. Fig. \ref{fig:fourier} indicates that for both $\rho$ and $T$, the numerical solution matches well with the reference solution for all Knudsen numbers. 

\begin{figure}[!hptb]
\centering
\subfigure[$r^{(n)}$]{
    \label{fig:Fourier_iter_Kndot01_a}
            \includegraphics[width = 0.23\textwidth]{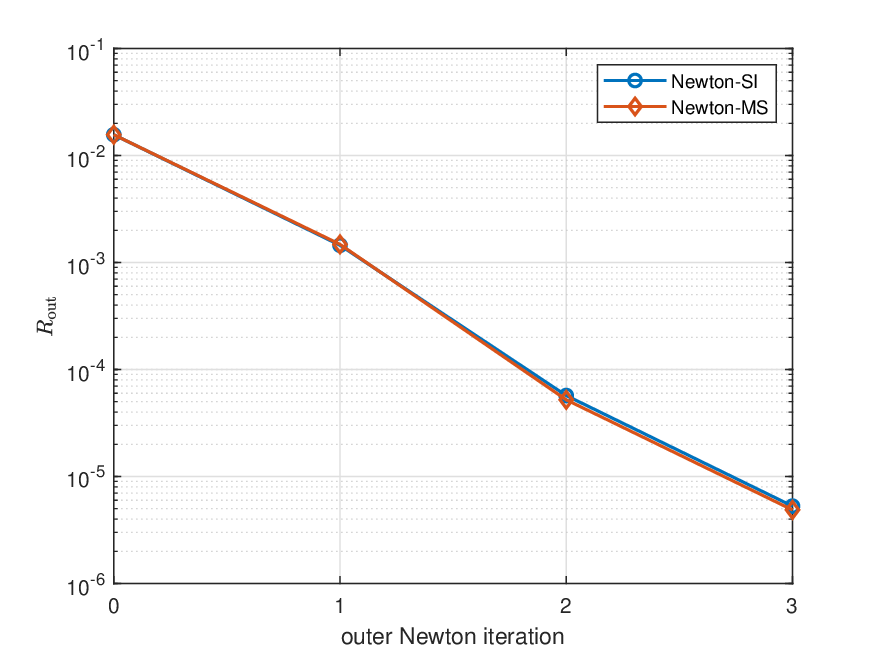}
}\hfill
\subfigure[$\mathcal{R}_{\rm in}^{(1,l)}$]{
    \label{fig:Fourier_iter_Kndot01_b}
            \includegraphics[width = 0.23\textwidth]{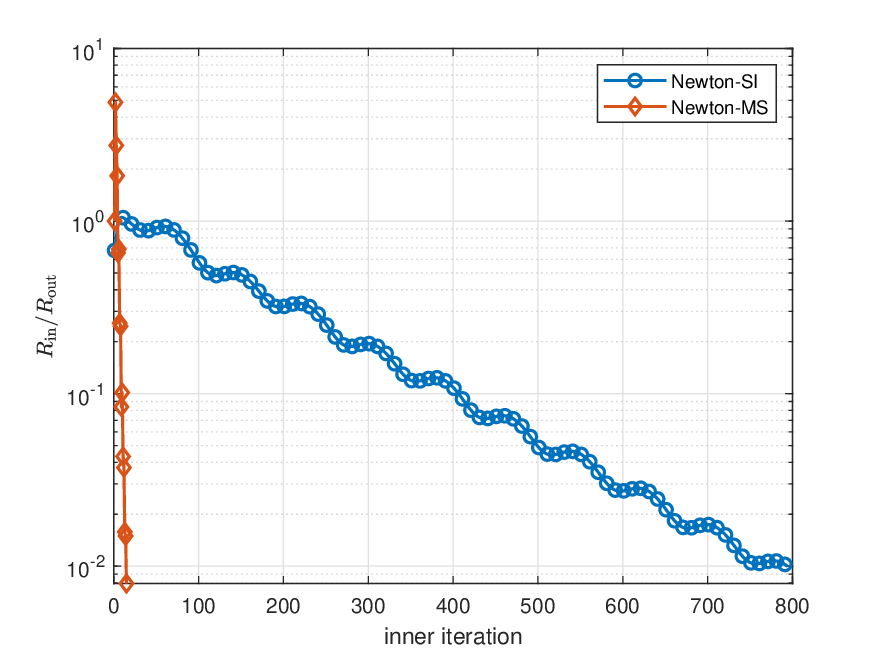}
}\hfill
\subfigure[$\mathcal{R}_{\rm in}^{(2,l)}$]{
    \label{fig:Fourier_iter_Kndot01_c}
            \includegraphics[width = 0.23\textwidth]{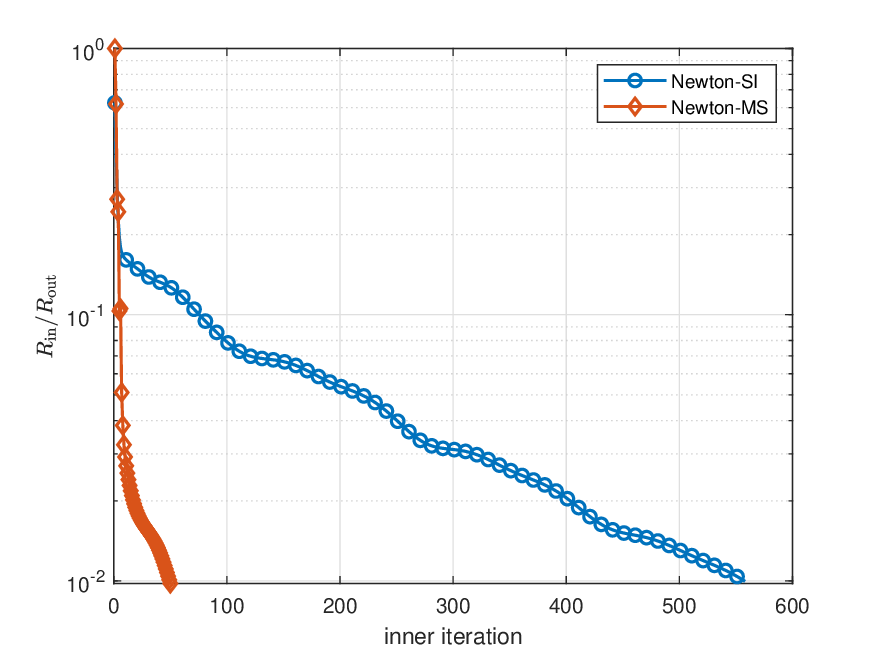}
}\hfill
\subfigure[$\mathcal{R}_{\rm in}^{(3,l)}$]{
    \label{fig:Fourier_iter_Kndot01_d}
            \includegraphics[width = 0.23\textwidth]{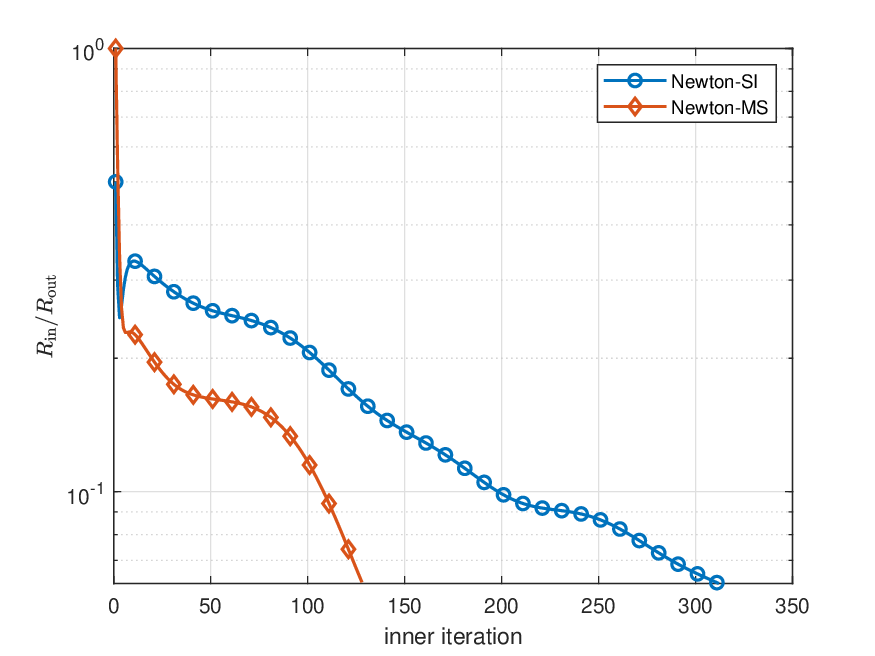}
}
\caption{(1D Fourier flow problem in Sec. \ref{sec:fourier}) Comparison of the convergence histories of outer Newton and inner iterations for Newton-MS and Newton-SI at $\Kn = 0.01$. (a) evolution of the outer Newton iteration residual $r^{(n)}$. (b) evolution of the inner iteration residual for the first outer Newton iteration $\mathcal{R}_{\rm in}^{(1,l)}$. (c) evolution of the inner iteration residual for the second outer Newton iteration $\mathcal{R}_{\rm in}^{(2,l)}$. (d) evolution of the inner iteration residual for the third outer Newton iteration $\mathcal{R}_{\rm in}^{(3,l)}$.}
    \label{fig:Fourier_iter_Kndot01}
\end{figure}

The efficiency comparison of Newton-MS, Newton-SI, and GSIS is summarized in Tab. \ref{tab:fourier_efficiency}. It indicates that for Newton-MS, the inner iteration number $N_{\rm in}$ remains small with decreasing $\Kn$, while for Newton-SI, $N_{\rm in}$ is increasing rapidly, and it even fails to converge when $\Kn = 0.001$. The behavior of $T_{\rm in}$ is similar, which is at the same order for $\Kn = 1$, but that of Newton-SI is $11$ times that of Newton-MS for $\Kn = 0.01$. Moreover, the total computational time $T_{\rm tol}$ of Newton-MS is also greatly reduced to more than $7$ times compared to Newton-SI for $\Kn = 0.01$. We can expect that efficiency can be further improved for smaller Knudsen numbers. Compared to GSIS, the outer iteration number is much smaller for Newton-MS. Thus, the total CPU time for Newton-MS is only half that of GSIS for $\Kn = 1$ and $0.001$, and about three-quarters for $\Kn = 0.1$ and $0.01$. This means that though the inner iteration is added, the total computational cost can still be reduced for Newton-MS compared to GSIS for all $\Kn$. Besides, Tab. \ref{tab:fourier_efficiency} also shows that for all Knudsen numbers, the computational cost for solving the macroscopic synthetic system is quite small, all less than $5\%$ of the total computational cost.

\begin{table}[htbp]
    \centering    
    \def\arraystrech{1.5}
    {\footnotesize
    \begin{tabular}{r l r r r r r r}
    \toprule
    \multirow{2}{*}{$\Kn$} & \multirow{2}{*}{method} & \multirow{2}{*}{$N_{\text{out}}$} & \multirow{2}{*}{$N_{\text{in}}$ } & \multicolumn{3}{c}{sub-time (s) (\%)} & \multirow{2}{*}{$T_{\rm tol}$} \\
    \cmidrule(lr){5-7}
          & & & & \makecell{$T_{\text{out}}$} & \makecell{$T_{\text{in}}$} & \makecell{$T_{\text{m}}$} & \\
    \midrule
    \multirow{3}{*}{1}
        & NMS  & 2 & 5  & 412.5 (78.8\%)  & 51.2 (9.80\%)    & 12.2 (2.3\%)   & 523.6  \\
          & NSI    & 2 & 7  & 410.6 (79.4\%)  & 66.3 (12.8\%)   & \NA          & 517.1  \\
          & GSIS         & 7 & - & 1037.2 (97.1\%) & \NA            & 3.1 (0.3\%)    & 1067.7 \\
    \midrule
    \multirow{3}{*}{0.1}
          & NMS  & 3 & 16 & 541.2 (49.0\%)  & 330.0 (29.9\%)  & 35.6 (3.2\%)   & 1103.5 \\
          & NSI    & 3 & 30 & 544.6 (44.6\%)  & 430.3 (35.2\%)  & \NA          & 1222.4 \\
          & GSIS         & 10 & -& 1479.3 (98.5\%) & \NA            & 4.3 (0.3\%)    & 1501.3 \\
    \midrule
    \multirow{3}{*}{0.01}
          & NMS  & 3 & 52  & 548.5 (27.5\%)  & 822.4 (41.2\%)  & 84.3 (4.2\%)   & 1996.9  \\
          & NSI    & 3 & 563 & 565.6 (3.60\%)   & 9428.9 (61.7\%) & \NA          & 15267.5 \\
          & GSIS         & 18 & -& 2516.8 (98.6\%) & \NA            & 7.0 (0.3\%)    & 2551.7  \\
    \midrule
    \multirow{3}{*}{0.001}
          & NMS  & 5 & 37  & 784.3 (36.0\%)  & 732.5 (33.7\%)  & 101.2 (4.6\%)  & 2176.5  \\
          & NSI    & -& -& \NA           & \NA            & \NA          & \NA   \\
          & GSIS         & 38 & - & 5162.8 (98.8\%) & \NA            & 18.6 (0.4\%)   & 5227.6  \\
    \bottomrule
    \end{tabular}
    }
    \caption{(1D Fourier flow problem in Sec. \ref{sec:fourier}) The efficiency comparison of Newton-MS (NMS), Newton-SI (NSI) and GSIS for different Knudsen numbers. For NMS and NSI, $N_{\mathrm{out}}$ and $N_{\mathrm{in}}$ denote the number of outer Newton iteration and the average number of inner iterations. For GSIS, $N_{\mathrm{out}}$ denotes the total number of iterations, and there is no $N_{\mathrm{in}}$. $T_{\text{out}}$ and $T_{\text{in}}$ are the computational time for obtaining $r^{(n)}$ \eqref{eq:r} and $\mathcal{R}_{\rm in}^{(n, l)}$ \eqref{eq:inner_residual_def}. $T_{\rm m}$ and $T_{\rm tol}$ are the computational time for solving macroscopic synthetic system \eqref{eq:dis_mac_eq3} and the total CPU time. }
    \label{tab:fourier_efficiency}
\end{table}

The convergence histories of the outer Newton and inner iterations for $\Kn = 0.01$ are plotted in Fig. \ref{fig:Fourier_iter_Kndot01}. It shows that the evolution of the outer Newton iteration residual $r^{(n)}$ for Newton-MS and Newton-SI is almost the same as in Fig. \ref{fig:Fourier_iter_Kndot01_a}, indicating that the macroscopic synthetic system does not deteriorate the convergence of the outer Newton iteration. Fig. \ref{fig:Fourier_iter_Kndot01_b} to \ref{fig:Fourier_iter_Kndot01_d} present the evolution of the inner iteration residual $\mathcal{R}_{\rm in}^{(n,l)}$ for each outer Newton iteration of both methods. It shows that compared to Newton-SI, the inner iteration residual of Newton-MS decreases much more quickly for each outer Newton iteration, which is also consistent with the results in Tab. \ref{tab:fourier_efficiency}.

\subsection{1D Couette flow problem}
\label{sec:couette}
In this section, the planner Couette flow problem is studied to validate the performance of Newton-MS in the shear-dominated problems. The scenario consists of two parallel plates at $x=0$ and $x=1$, which move in opposite tangential directions with wall velocities $u_w=\pm 0.5$ and fixed temperature $T_w=1.0$. Fully diffuse reflection boundary conditions are imposed on both walls. The similar example is also tested in \cite{jaiswal2019discontinuous}. The discretization of the spatial and microscopic velocity space is the same as that in Sec. \ref{sec:fourier}.  

\begin{figure}[!hptb]
\centering
\subfigure[$\rho$]{
            \includegraphics[width = 0.3\textwidth]{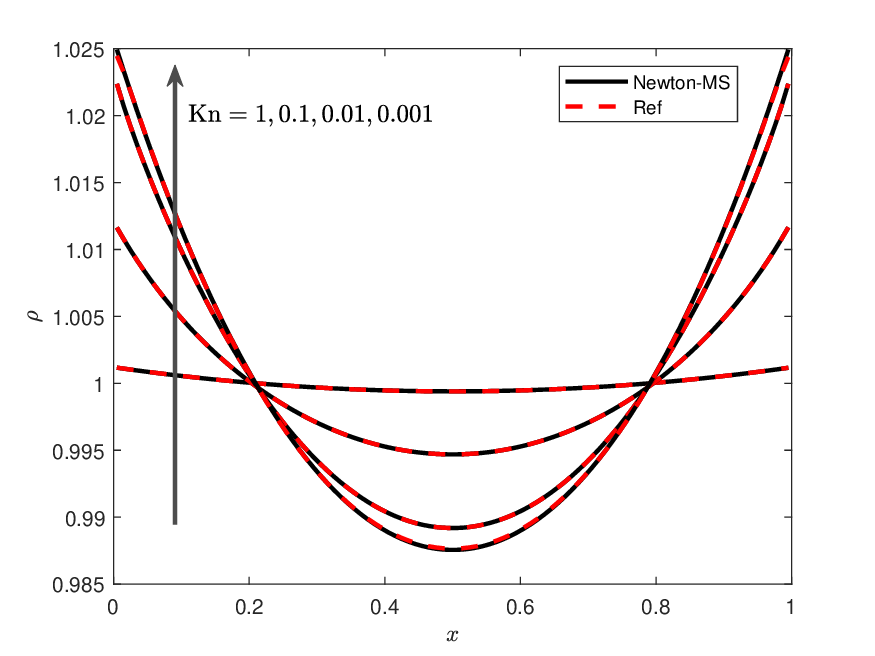}
}
\subfigure[$u_y$]{
            \includegraphics[width = 0.3\textwidth]{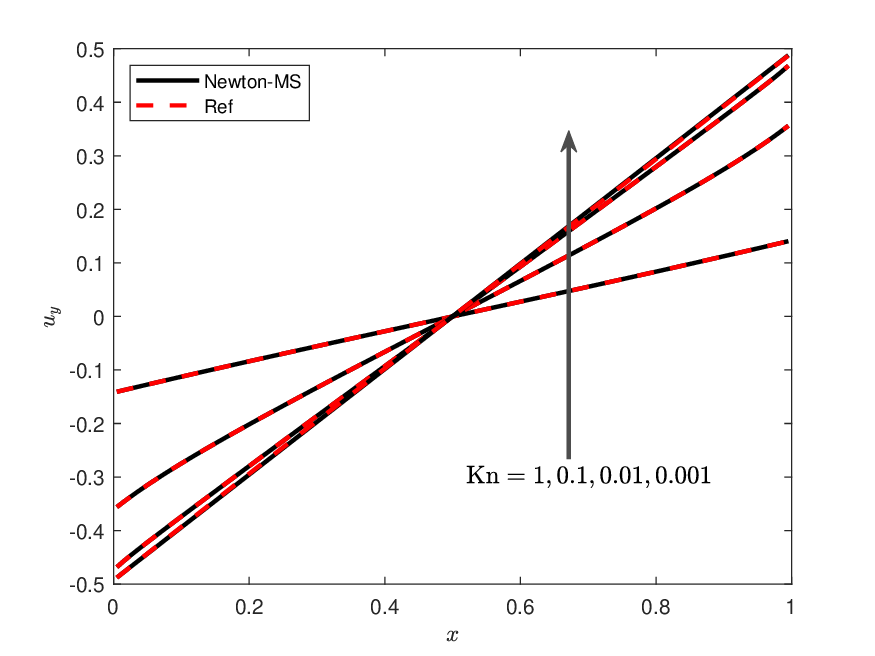}
}
\subfigure[$T$]{
            \includegraphics[width = 0.3\textwidth]{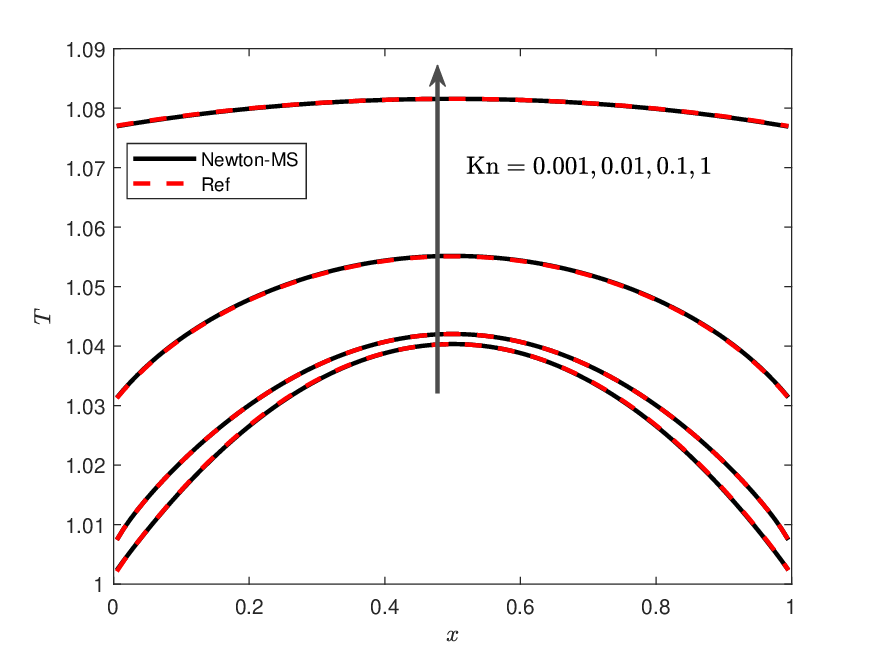}
}
\caption{(1D Couette flow in Sec. \ref{sec:couette}) Steady state numerical solution of the density $\rho$, macroscopic velocity in $y$-axes $u_y$, and temperature $T$ for $\Kn=1,0.1,0.01$, and $0.001$. Here, the solid black lines are the numerical solution by Newton-MS, and the red dashed lines are those of the reference solution.}
\label{fig:couette}
\end{figure}

The steady state numerical solution of the density $\rho$, the tangential velocity $u_y$ and the temperature for $\Kn = 1, 0.1, 0.01$ and $0.001$ is shown in Fig. \ref{fig:couette}, where the reference solution is also plotted. Here,  the reference solution of  $\Kn = 0.001$ is obtained by the compressible Navier-Stokes equations with slip boundary conditions, while the others are obtained by Newton-SI. Fig. \ref{fig:couette} shows that for all Knudsen numbers, the numerical solution agrees well with the reference solution. 

\begin{table}[htbp]
    \centering    
    \def\arraystrech{1.5}
    {\footnotesize
    \begin{tabular}{r l r r r r r r}
    \toprule
    \multirow{2}{*}{$\Kn$} & \multirow{2}{*}{method} 
    & \multirow{2}{*}{$N_{\text{out}}$} 
    & \multirow{2}{*}{\makecell{$N_{\text{in}}$}} 
    & \multicolumn{3}{c}{sub-time (s) (\%)} 
    & \multirow{2}{*}{\makecell{$T_{\text{tol}}$}} \\
    \cmidrule(lr){5-7}
          & & & & \makecell{$T_{\text{out}}$} & \makecell{$T_{\text{in}}$} & \makecell{$T_{\text{m}}$} & \\
    \midrule
    \multirow{3}{*}{1}
          & NMS  & 2 & 6    & 416.4 (78.4\%)  & 65.3 (12.3\%)   & 11.9 (2.2\%)   & 531.1  \\
          & NSI    & 2 & 6    & 419.5 (80.4\%)  & 64.1 (12.3\%)   & \NA          & 521.6  \\
          & GSIS         & 8 & {-} & 1221.3 (97.0\%) & \NA            & 3.7 (0.3\%)    & 1258.9 \\
    \midrule
    \multirow{3}{*}{0.1}
          & NMS  & 2  & 11   & 413.8 (69.9\%)  & 120.7 (20.4\%)  & 9.7 (1.6\%)    & 591.9  \\
          & NSI    & 2  & 26   & 421.6 (49.4\%)  & 271.6 (31.8\%)  & \NA          & 853.4  \\
          & GSIS         & 14 & {-} & 1909.6 (94.8\%) & \NA            & 7.6 (0.4\%)    & 2013.9 \\
    \midrule
    \multirow{3}{*}{0.01}
          & NMS  & 2  & 38   & 421.6 (38.3\%)  & 412.7 (37.5\%)  & 33.9(3.1\%)   & 1101.3  \\
          & NSI    & 2  & 426  & 431.2 (5.5\%)   & 4207.6 (54.5\%) & \NA          & 7721.6 \\
          & GSIS         & 23 & {-} & 3220.8 (97.9\%) & \NA            & 12.9 (0.4\%)   & 3289.7 \\
    \midrule
    \multirow{3}{*}{0.001}
          & NMS  & 4  & 31   & 712.4 (34.3\%)  & 811.3 (39.0\%)  & 91.6 (4.4\%)   & 2079.3 \\
            & NSI    & {-} & {-} & \NA           & \NA            & \NA          & \NA   \\
          & GSIS         & 56 & {-} & 7386.4 (98.8\%) & \NA            & 29.6 (0.4\%)   & 7478.1 \\
    \bottomrule
    \end{tabular}
    }
    \caption{(1D Couette flow problem in Sec. \ref{sec:couette}) 
     The efficiency comparison of Newton-MS (NMS), Newton-SI (NSI) and GSIS for different Knudsen numbers. For NMS and NSI, $N_{\mathrm{out}}$ and $N_{\mathrm{in}}$ denote the number of outer Newton iteration and the average number of inner iterations. For GSIS, $N_{\mathrm{out}}$ denotes the total number of iterations, and there is no $N_{\mathrm{in}}$. $T_{\text{out}}$ and $T_{\text{in}}$ are the computational time for obtaining $r^{(n)}$ \eqref{eq:r} and $\mathcal{R}_{\rm in}^{(n, l)}$ \eqref{eq:inner_residual_def}. $T_{\rm m}$ and $T_{\rm tol}$ are the computational time for solving macroscopic synthetic system \eqref{eq:dis_mac_eq3} and the total CPU time.}
    \label{tab:couette_efficiency}
\end{table}

The comparison of the computational efficiency between Newton-MS, Newton-SI, and GSIS is summarized in Tab. \ref{tab:couette_efficiency}, where the similar behavior is found as in Sec. \ref{sec:fourier}. For $\Kn=1$, Newton-MS has a comparable cost to Newton-SI, since the source iteration is already efficient in the rarefied regime and the additional macroscopic correction brings only limited benefit. As $\Kn$ decreases, however, the advantage of Newton-MS becomes increasingly evident. Newton-SI suffers from the slow convergence of the inner source iteration, whereas Newton-MS keeps the number of inner iterations small.  Compared with GSIS, Newton-MS also requires fewer outer iterations and achieves a lower total CPU time, especially in the near-continuum regime. The computational cost $T_{\rm m}$ to solve the macroscopic synthetic system all remains a small fraction of the total cost for all Knudsen numbers.

\begin{figure}[!hptb]
\centering
\subfigure[$r^{(n)}$]{
\label{fig:couette_iter_Kndot01_a}
    \includegraphics[width = 0.3\textwidth]{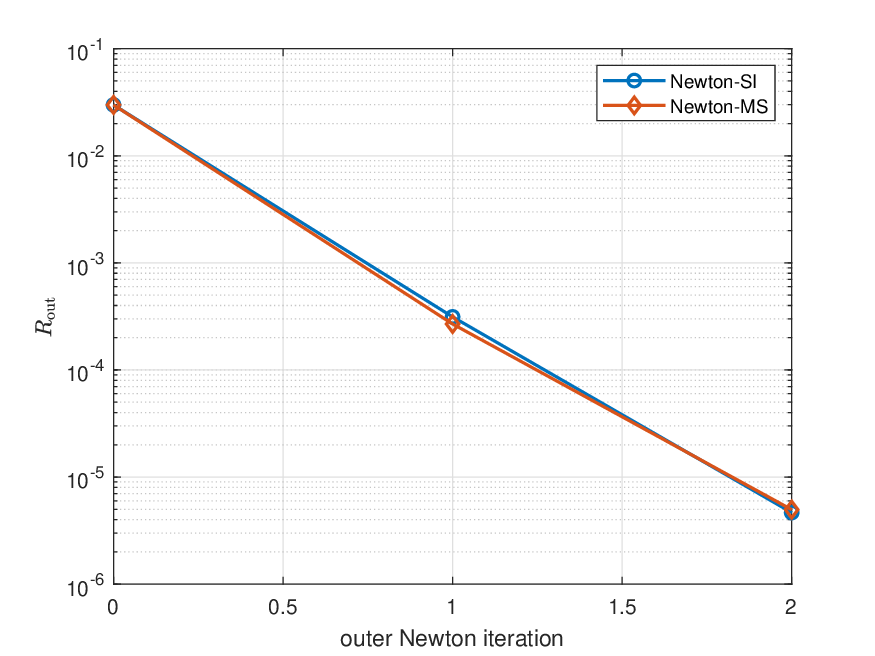}
}
\subfigure[$\mathcal{R}_{\rm in}^{(1,l)}$]{
\label{fig:couette_iter_Kndot01_b}
            \includegraphics[width = 0.3\textwidth]{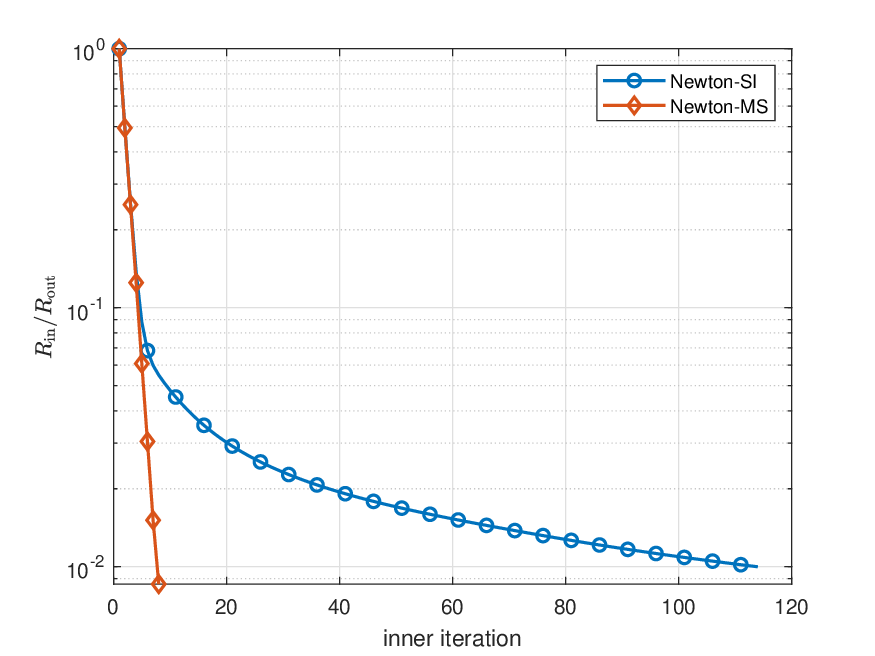}
}
\subfigure[$\mathcal{R}_{\rm in}^{(2,l)}$]{
\label{fig:couette_iter_Kndot01_c}
            \includegraphics[width = 0.3\textwidth]{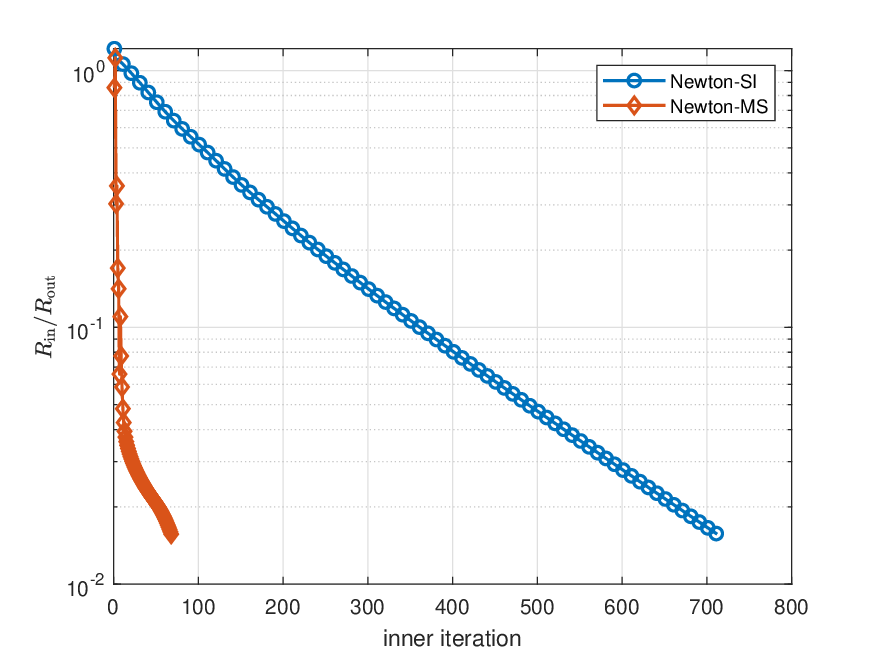}
}
\caption{(1D Couette flow problem in Sec. \ref{sec:couette}) Comparison of the convergence histories of outer Newton and inner iterations for Newton-MS and Newton-SI at $\Kn = 0.01$. (a) evolution of the outer Newton iteration residual $r^{(n)}$. (b) evolution of the inner iteration residual for the first outer Newton iteration $\mathcal{R}_{\rm in}^{(1,l)}$. (c) evolution of the inner iteration residual for the second outer Newton iteration $\mathcal{R}_{\rm in}^{(2,l)}$. }
\label{fig:couette_iter_Kndot01}
\end{figure}

The convergence histories of the outer Newton and inner iterations for $\Kn = 0.01$ are plotted in Fig. \ref{fig:couette_iter_Kndot01}, the behavior of which is similar to that in Sec. \ref{sec:fourier}. Precisely, the evolution of the outer Newton iteration residual $r^{(n)}$ for Newton-MS and Newton-SI is almost the same as in Fig. \ref{fig:couette_iter_Kndot01_a}, indicating that the macroscopic synthetic system does not deteriorate the convergence of the outer Newton iteration. Fig. \ref{fig:couette_iter_Kndot01_b} to \ref{fig:couette_iter_Kndot01_c} present the evolution of the inner iteration residual $\mathcal{R}_{\rm in}^{(n,l)}$ for each outer Newton iteration of both methods. It shows that compared to Newton-SI, the inner iteration residual of Newton-MS decreases much more quickly for each outer Newton iteration, which is also consistent with the results in Tab. \ref{tab:couette_efficiency}. The inner iteration residual reaches the tolerance much faster for Newton-MS in each outer Newton iteration, compared to Newton-SI, indicating the effect of the acceleration for the macroscopic synthetic system when applied to the shear-dominated problems.

\subsection{2D lid-driven cavity flow problem}
\label{sec:cavity}
\begin{figure}[hptb]
\centering
\subfigure[density $\rho$]{
        \includegraphics[width = 0.4\textwidth]{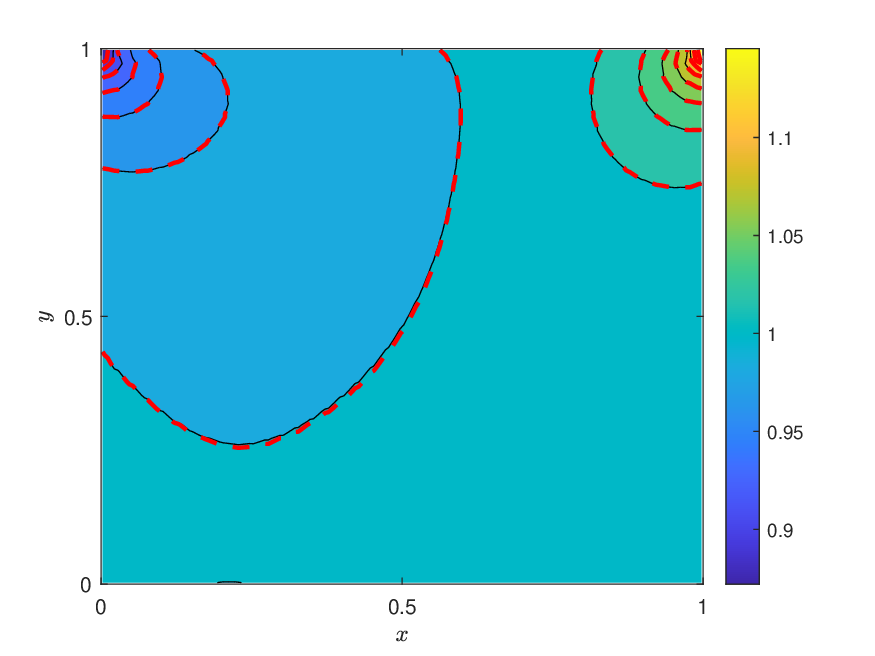}
}\qquad 
\subfigure[temperature $T$]{
        \includegraphics[width = 0.4\textwidth]{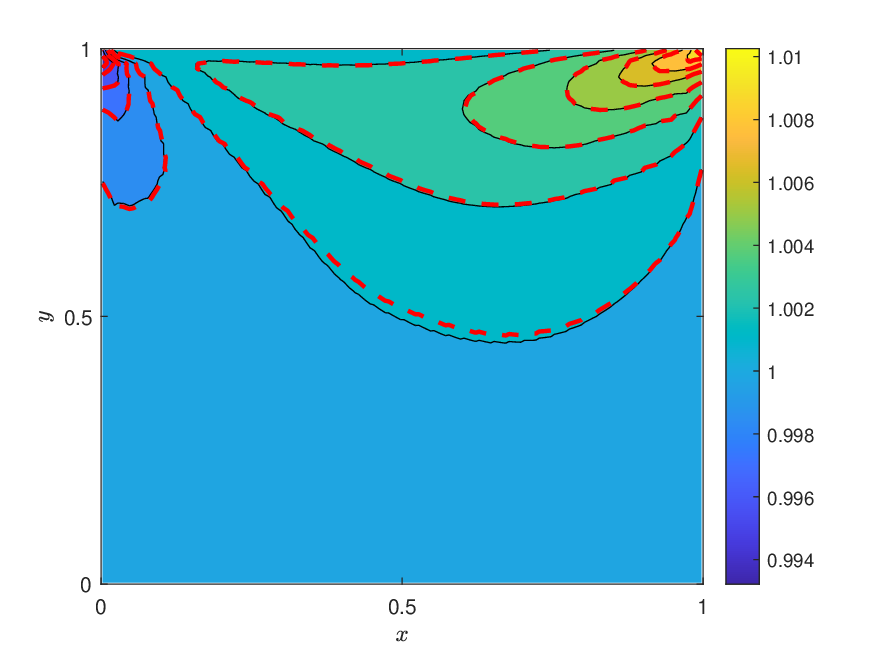}
}\\
\subfigure[velocity in $x$-direction $u_x$]{
        \includegraphics[width = 0.4\textwidth]{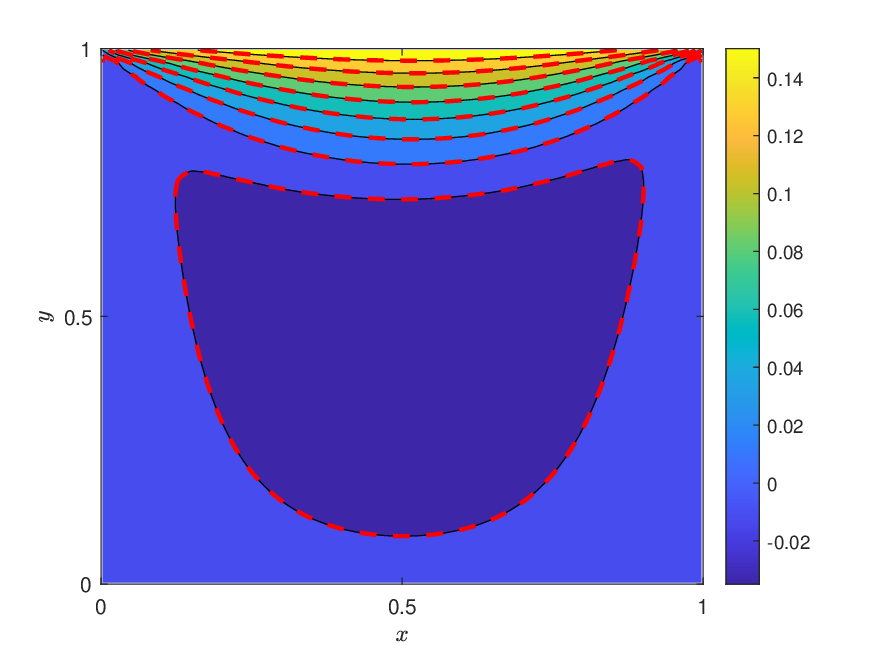}
}\qquad 
\subfigure[velocity in $y$-direction $u_y$]{
        \includegraphics[width = 0.4\textwidth]{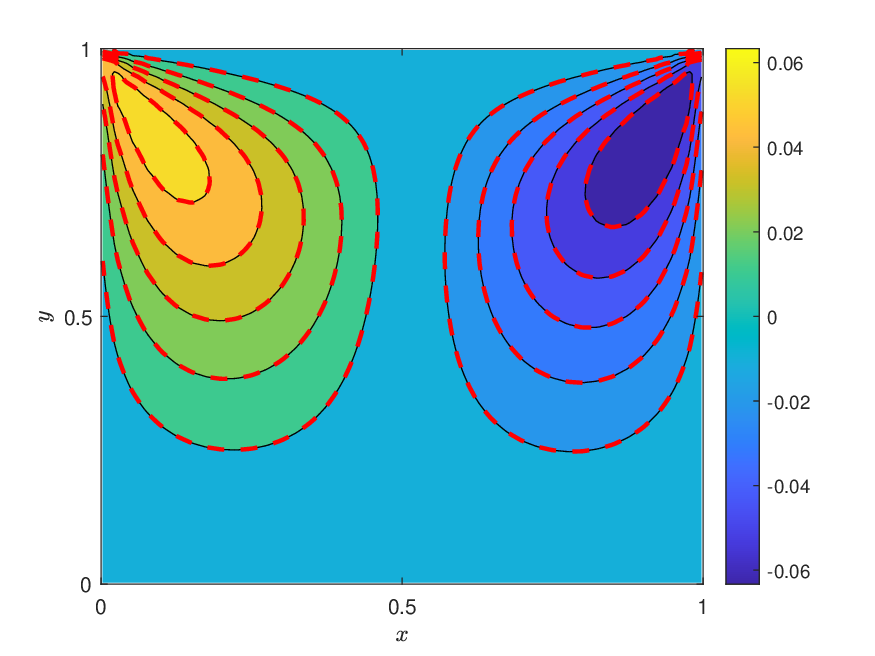}
}
\caption{(2D lid-driven cavity flow problem in Sec. \ref{sec:cavity}) Steady state numerical solution of the density $\rho$, temperature $T$, and macroscopic velocity $u_x$ and $u_y$ for $\Kn=0.01$. Here, the solid black lines are the numerical solution by Newton-MS, and the red dashed lines are those of the reference solution.}
\label{fig:Cavity}
\end{figure}
In this section, the two-dimensional lid-driven cavity flow problem is studied to demonstrate the capability of Newton-MS for the high dimensional problems. The scenario is a unit square domain $(x, y) \in [0,1]^2$. The top wall at $y = 1$ is moving with a constant tangential velocity $\bu_w = (0.5, 0, 0)$, while the other three walls are stationary. All walls are maintained at the fixed temperature $T_w=1.0$, and fully diffuse reflection boundary conditions are imposed. Similar problems are also studied in \cite{naris2005driven,jaiswal2019discontinuous}. The mesh size in the spatial space is $N_{\rm el} = 40 \times 40$ with the DG polynomial degree set to $N_p = 2$. For the microscopic velocity space, the number of Fourier modes is chosen as $N = 16$.

\begin{figure}[!hptb]
\centering
\subfigure[$r^{(n)}$]{
\label{fig:cavity_iter_a}
            \includegraphics[width = 0.23\textwidth]{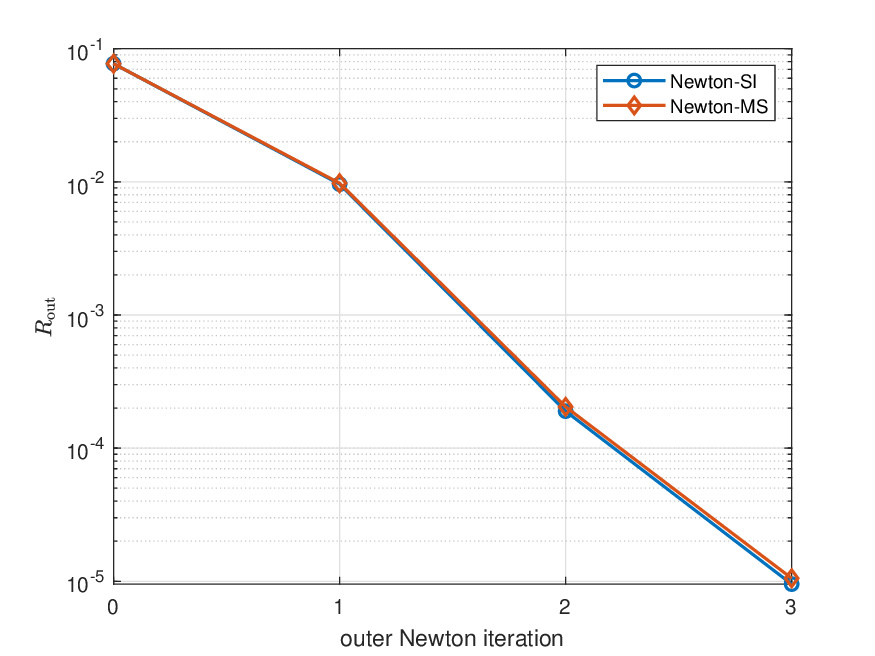}
}
\subfigure[$\mathcal{R}_{\rm in}^{(1,l)}$]{
\label{fig:cavity_iter_b}
            \includegraphics[width = 0.23\textwidth]{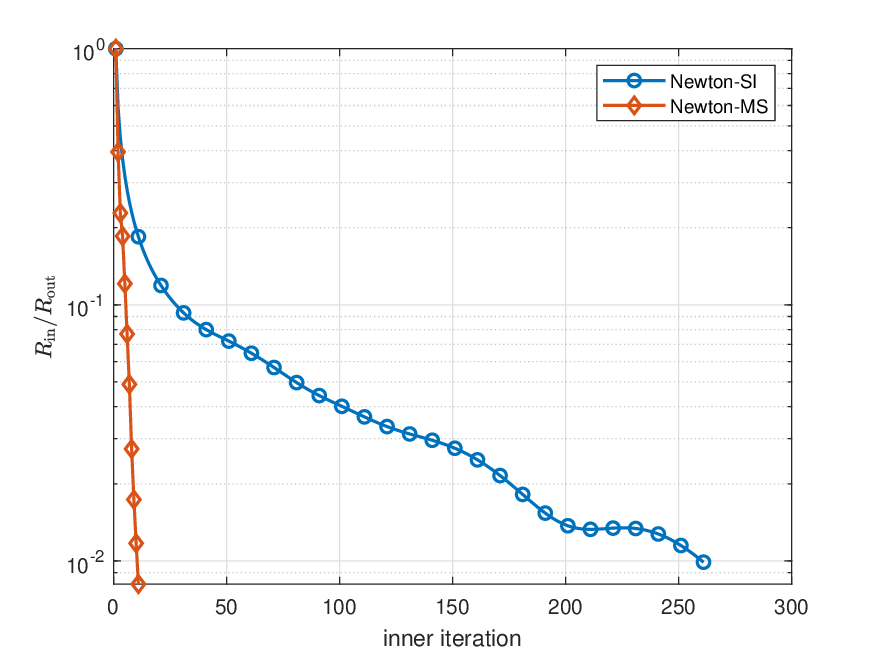}
}
\subfigure[$\mathcal{R}_{\rm in}^{(2,l)}$]{
\label{fig:cavity_iter_c}
            \includegraphics[width = 0.23\textwidth]{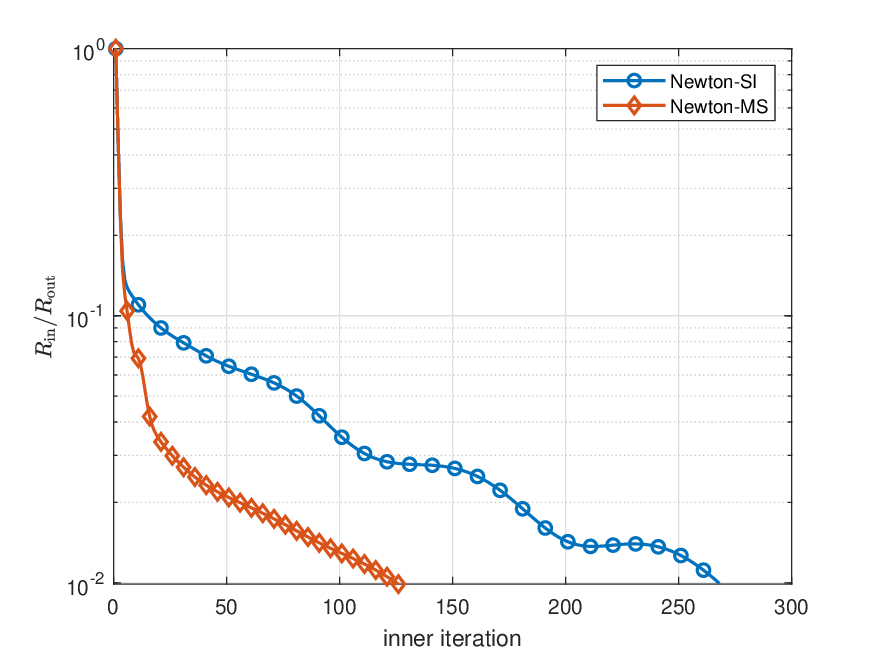}
}
\subfigure[$\mathcal{R}_{\rm in}^{(3,l)}$]{
\label{fig:cavity_iter_d}
            \includegraphics[width = 0.23\textwidth]{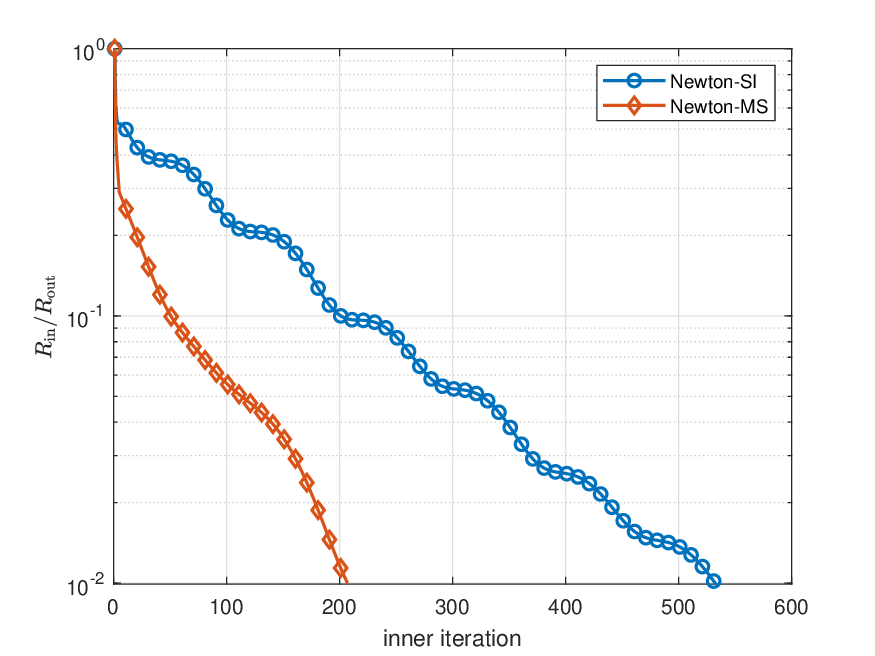}
}
\caption{(2D lid-driven cavity flow problem in Sec. \ref{sec:cavity}) Comparison of the convergence histories of outer Newton and inner iterations for Newton-MS and Newton-SI at $\Kn = 0.01$.
}
\label{fig:cavity_iter}
\end{figure}

The steady state numerical solution of the density $\rho$, temperature $T$, the macroscopic velocity $u_x$ and $u_y$ for $\Kn = 0.01$ is shown in Fig. \ref{fig:Cavity}, where the reference solution obtained by Newton-SI is also plotted. It shows that the numerical solution matches well with the reference solution. The computational efficiency of Newton-MS and Newton-SI is summarized in Tab. \ref{tab:cavity_perf}. It shows that though the number of outer Newton iteration is the same for both methods, the inner iteration number of Newton-MS is much smaller for each outer Newton iteration. Therefore, the total computational time of Newton-MS is less than half of Newton-SI. Moreover, the additional cost of solving the macroscopic synthetic system is quite small, less than $5\%$ of the total cost. 
\begin{table}[htbp]
    \centering
 \def\arraystrech{1.5}
     {\footnotesize
     \begin{tabular}{c l c r c c c r c c}
    \toprule
    \multirow{2}{*}{$\Kn$} & \multirow{2}{*}{method} & \multirow{2}{*}{$N_{\text{out}}$} & \multicolumn{3}{c}{$N_{\text{in}}$} &\multicolumn{3}{c}{sub-time (h)}& \multirow{2}{*}{\makecell{$T_{\text{tol}}$ (h)}} \\
    \cmidrule(lr){4-6}
    \cmidrule(lr){7-9}
          & & & \multicolumn{1}{c}{1} & \multicolumn{1}{c}{2} & \multicolumn{1}{c}{3}& \multicolumn{1}{c}{\makecell{$T_{\text{out}}$}} & \multicolumn{1}{c}{\makecell{$T_{\text{in}}$}} & \multicolumn{1}{c}{\makecell{$T_{\text{m}}$}} & \\
    \midrule
    \multirow{2}{*}{0.01}
        & NMS   & 3   & 10   & 125 & 218 &2.3& 9.6 &0.9 &18.5 \\
       & NSI    & 3  & 262 & 268 &621 & 2.3 &28.5&-&46.8\\
    \bottomrule
    \end{tabular}
    }
      \caption{(2D lid-driven cavity flow problem) The efficiency comparison of Newton-MS (NMS), and Newton-SI (NSI) for $\Kn = 0.01$.}
       \label{tab:cavity_perf}
\end{table}
The convergence histories of the outer Newton and inner iterations for Kn = 0.01 are plotted in Fig. \ref{fig:cavity_iter}. The evolution of the outer Newton residuals for Newton-MS and Newton-SI is similar, as shown in Fig. \ref{fig:cavity_iter_a}. The evolution of $\mathcal{R}_{\rm in}^{(n,l)}$ for each Newton iteration are presented in Fig. \ref{fig:cavity_iter_b} to \ref{fig:cavity_iter_d}, respectively. They illustrate that the decay of $\mathcal{R}_{\rm in}^{(n,l)}$ is much faster for Newton-MS compared to Newton-SI, which is also consistent with the results in Tab. \ref{tab:cavity_perf}.




\subsection{2D thermal cavity flow problem}
\label{sec:tcavity}
In this section, the 2D thermal cavity flow problem is studied. The scenario is also a unit square domain $(x,y)\in[0,1]^2$ as in Sec. \ref{sec:cavity}, but all four walls are stationary. The temperature of the top wall is set as $T_{\mathrm{top}}=1.2$, while all others are fixed as $T = 1$. The fully diffuse reflection boundary conditions are imposed on all walls. Different from the lid-driven cavity flow in Sec. \ref{sec:cavity}, there is no external momentum added, and the flow is driven entirely by the difference of the wall temperature. The similar example is also studied in \cite{jaiswal2019discontinuous}. 

\begin{figure}[!hptb]
\centering
\subfigure[$\rho$]{
            \includegraphics[width = 0.45\textwidth]{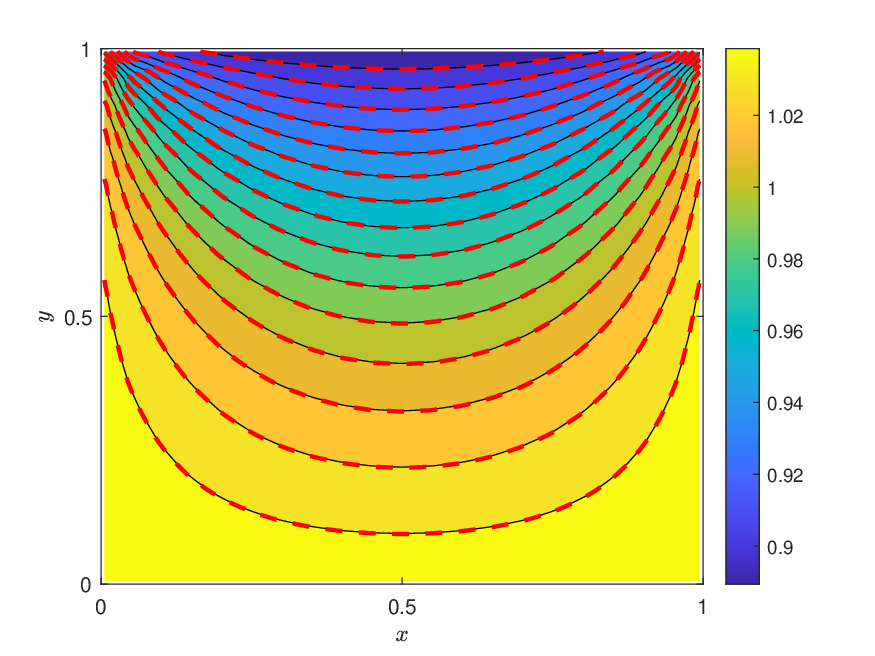}
}
\subfigure[$T$]{
            \includegraphics[width = 0.45\textwidth]{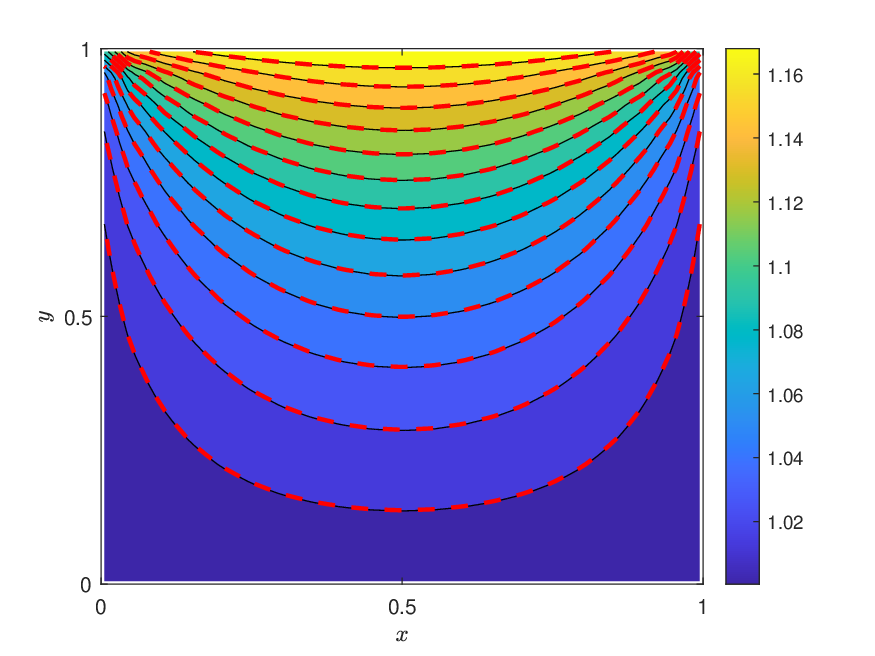}
}
\caption{(2D thermal cavity flow problem in Sec. \ref{sec:tcavity}) Steady state numerical solution of the density $\rho$, and temperature $T$ for $\Kn=0.01$. Here, the solid black lines are the numerical solution by Newton-MS, and the red dashed lines are those of the reference solution.}
\label{fig:Tcavity}
\end{figure}

The same discretization of the spatial space and microscopic velocity as in Sec. \ref{sec:cavity} is adopted. The steady state numerical solution of the density $\rho$, and the temperature $T$ for $\Kn = 0.01$ is plotted in Fig. \ref{fig:Tcavity}, where it matches well with the reference solution obtained by Newton-SI. The computational efficiency is summarized in Tab.~\ref{tab:tcavity_perf}, where the behavior of both methods is similar as that in Sec. \ref{sec:cavity}. The total computational time of Newton-MS is reduced to one third of that by Newton-SI, which is also mainly due to the reduced number of the inner iterations. Moreover, the additional cost of solving the macroscopic synthetic system still remains quite small, compared to the total computational time saved.


\begin{figure}[!hptb]
\centering
\subfigure[$r^{(n)}$]{
\label{fig:Tcavity_iter_a}
            \includegraphics[width = 0.18\textwidth]{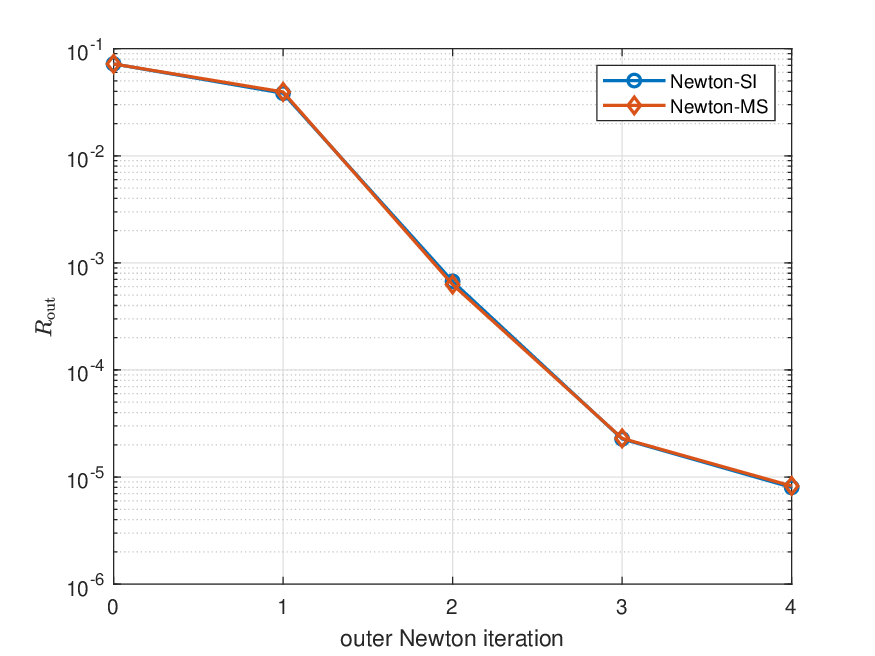}
}
\subfigure[$\mathcal{R}_{\rm in}^{(1,l)}$]{
\label{fig:Tcavity_iter_b}
            \includegraphics[width = 0.18\textwidth]{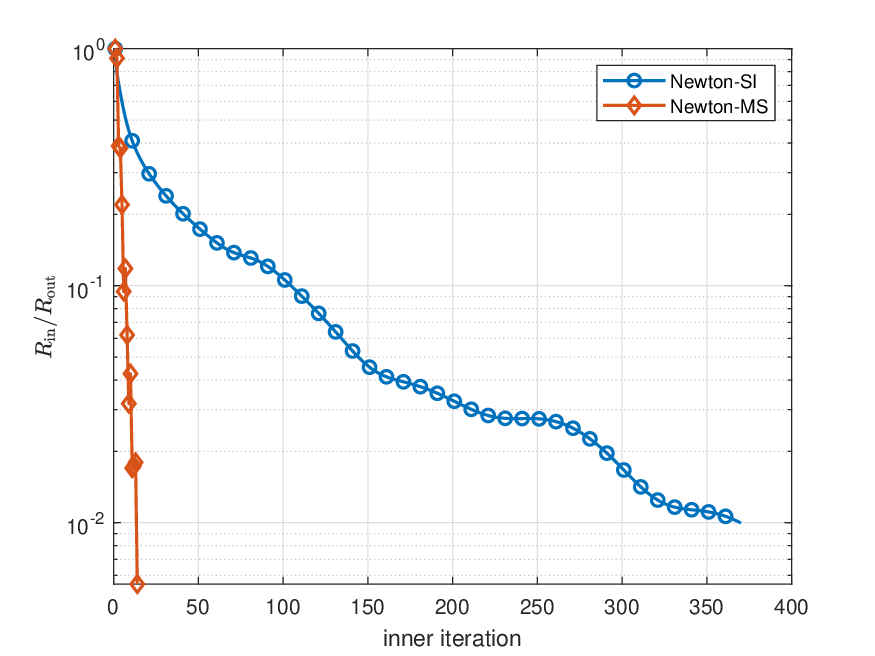}
}
\subfigure[$\mathcal{R}_{\rm in}^{(2,l)}$]{
\label{fig:Tcavity_iter_c}
            \includegraphics[width = 0.18\textwidth]{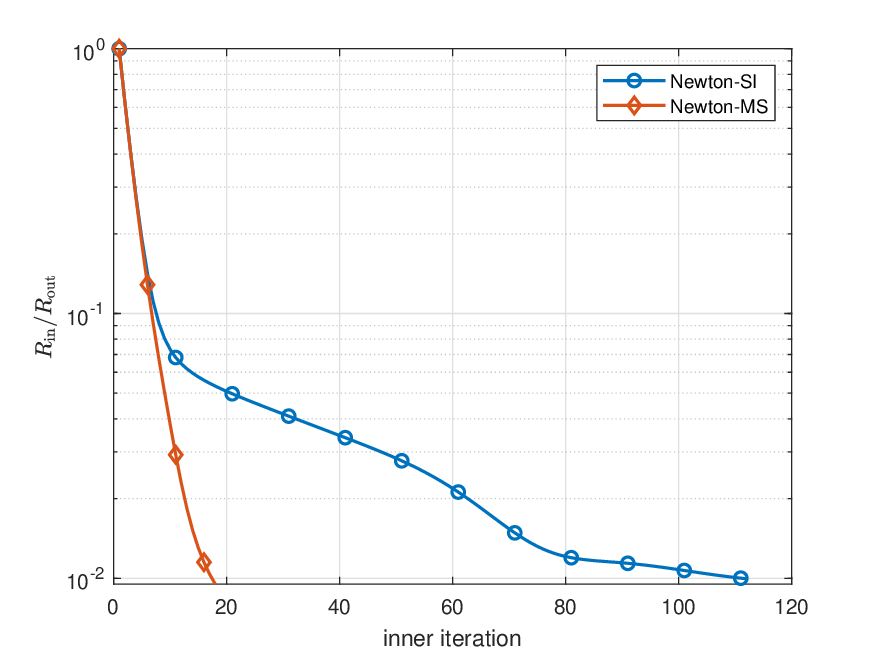}
}
\subfigure[$\mathcal{R}_{\rm in}^{(3,l)}$]{
\label{fig:Tcavity_iter_d}
            \includegraphics[width = 0.18\textwidth]{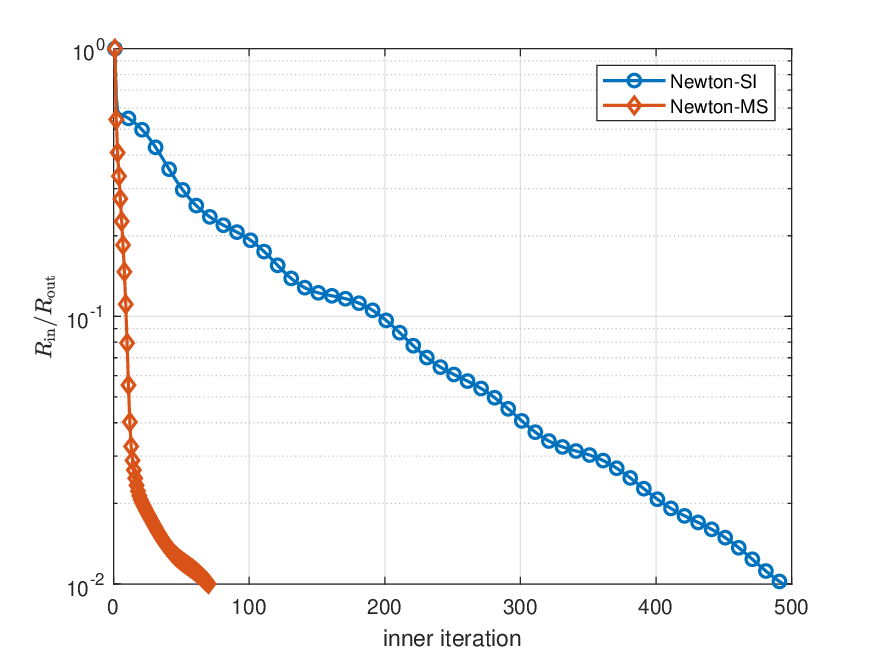}
}
\subfigure[$\mathcal{R}_{\rm in}^{(4,l)}$]{
\label{fig:Tcavity_iter_e}
            \includegraphics[width = 0.18\textwidth]{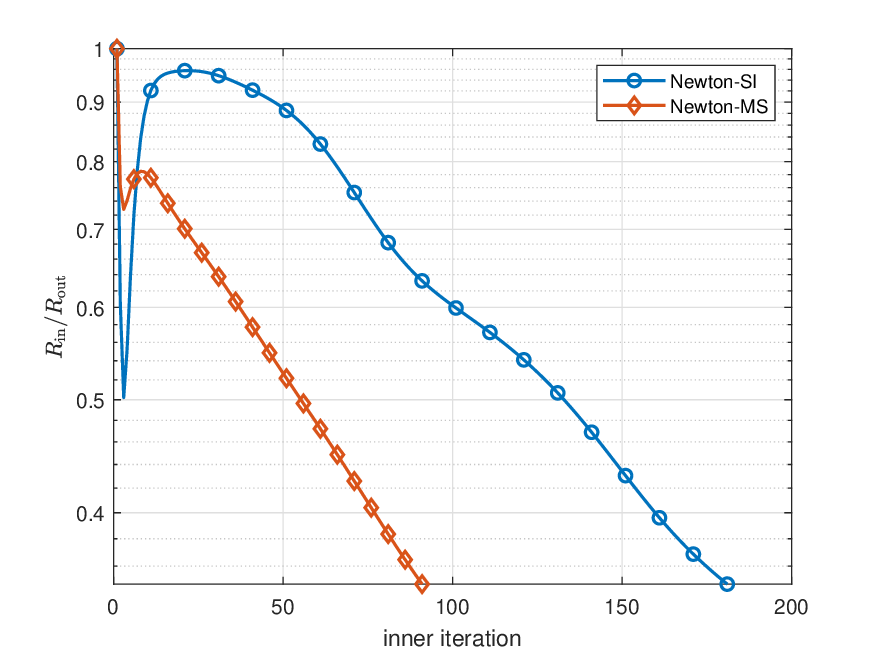}
}
\caption{(2D thermal cavity flow problem in Sec. \ref{sec:tcavity}) Comparison of the convergence histories of outer Newton and inner iterations for Newton-MS and Newton-SI at $\Kn = 0.01$. 
}
\label{fig:Tcavity_iter}
\end{figure}

The convergence histories of the outer Newton and inner iterations for $\Kn = 0.01$ are plotted in Fig. \ref{fig:Tcavity_iter}. The outer Newton residuals decay at almost the same rates for Newton-SI and Newton-MS, as shown in Fig. \ref{fig:Tcavity_iter_a}, indicating that the synthetic acceleration does not affect the outer nonlinear convergence. For the inner iteration, the evolution of $\mathcal{R}_{\rm in}^{(n,l)}$ is the similar to that in Sec. \ref{sec:cavity}, which all shows a much more rapid decay rate for Newton-MS compared to that of Newton-SI.


\begin{table}[htbp]
    \centering
     \def\arraystrech{1.5}
     {\footnotesize
    \begin{tabular}{c l c r r r r c r c c}
    \toprule
    \multirow{2}{*}{$\Kn$} & \multirow{2}{*}{method} & \multirow{2}{*}{$N_{\mathrm{out}}$} & \multicolumn{4}{c}{$N_{\mathrm{in}}$} & \multicolumn{3}{c}{sub-time (h)} & \multirow{2}{*}{\makecell{$T_{\text{tol}}$}} \\
    \cmidrule(lr){4-7}
    \cmidrule(lr){8-10}
          & & & \multicolumn{1}{c}{1} & \multicolumn{1}{c}{2} & \multicolumn{1}{c}{3} & \multicolumn{1}{c}{4} & \multicolumn{1}{c}{\makecell{$T_{\mathrm{out}}$}} & \multicolumn{1}{c}{\makecell{$T_{\mathrm{in}}$}} & \multicolumn{1}{c}{\makecell{$T_{\mathrm{m}}$}} & \\
    \midrule
    \multirow{2}{*}{0.01}
       & NMS & 4 & 14  & 18  & 70  & 91  & 2.9 & 6.6  & 0.7   & 16.2 \\
       & NSI   & 4 & 370 & 112 & 494 & 181 & 2.9 & 30.2 & -    & 49.3 \\
    \bottomrule
    \end{tabular}
    }
    \caption{(2D thermal cavity flow problem) The efficiency comparison of Newton-MS (NMS), and Newton-SI (NSI) for $\Kn = 0.01$.}
    \label{tab:tcavity_perf}
\end{table}
\section{Conclusion}
\label{sec:con}
In this work, a modified Newton's method accelerated by a macroscopic synthetic system (Newton-MS) is proposed for the steady-state Boltzmann equation. For Newton-MS, the outer iteration is the normal Newton method, while the linearized collision operator is utilized instead of the quadratic Boltzmann collision operator in the inner iteration to reduce computational cost. For the inner iteration, a macroscopic synthetic system based on a Chapman-Enskog closure is first derived, and then the information of the updated macroscopic variables obtained by solving this macroscopic synthetic system is added in the general source iteration to accelerate inner iteration convergence. The discrete matrix form of the macroscopic synthetic system is derived in the framework of the discontinuous Galerkin method for the practical implementation, where the computational cost can be greatly reduced compared to directly discretizing the continuous macroscopic synthetic system. Numerical results indicate that Newton-MS maintains consistent convergence properties from the free-molecular regime to the continuum limit and offers a practical and effective alternative to standard Newton-SI schemes for rarefied gas simulations.

\section*{Acknowledgments}
  We thank Prof. Lei Wu from Southern University of Science and Technology for providing the GSIS code. This work of Yanli Wang is partially supported by the Science Challenge program (NO. TZ2025016). This work of Pei Zhang was supported by the China Scholarship Council (NO. 202404890001). The work of Zhenning Cai was supported by the Academic Research Fund of the Ministry of Education of Singapore under grant A-8002392-00-00.
\appendix

\section{Derivation of \texorpdfstring{$\sigma^{(\nabla)}$}{sigma{(nabla)}} and \texorpdfstring{ $q^{(\nabla)}$}{q{(nabla)}}}
\label{sec:constitutive_derivation}
In this appendix, the detailed deduction of the stress tensor $\sigma^{(\nabla)}$ and heat flux $q^{(\nabla)}$ in \eqref{eq:sigma_q_g} is presented. The derivation is based on the expansion of $\bv\cdot\nabla_{\bx}\mathcal{E}$ on the microscopic velocity space, and the rotational invariance of the linear collision operator $\mL$.

\subsection{Deducing \texorpdfstring{$\sigma^{(\nabla)}$}{sigma{(nabla)}} and \texorpdfstring{ $q^{(\nabla)}$}{q{(nabla)}}}
We first introduce some abbreviations for convenience as 
\begin{subequations}
  \label{eq:grad_shorthand}
  \begin{gather}
       \label{eq:abbre}
  \partial_i = \pd{}{x_i}, \qquad   A_i = \frac{\partial_i\rho}{\rho},
  \qquad
  U_{ij} = \partial_i \bu_{j},
  \qquad
  G_i = \partial_i T,
  \qquad
  D_{ij} = \partial_i\delta \bu_j, \\
   \mathcal{A}_{ij}(\bc) = c_i c_j - \frac{|\bc|^2}{3}\delta_{ij},
  \qquad
  \mathcal{B}_i(\bc) = \frac{1}{2}|\bc|^2 c_i,
\qquad 
   \alpha = \frac{\delta\rho}{\rho},
  \qquad
  \beta_i = \frac{\delta u_i}{T},
  \qquad
  \vartheta = \frac{\delta T}{2T^2}.
  \label{eq:abg_def}
\end{gather}
\end{subequations}
Then the local equilibrium \eqref{eq:linear_Max} can be rewritten as 
\begin{equation}
  \label{eq:H_def} 
    \mathcal{E}(\bx,\bv) = \mathcal{M}\,\chi,
  \qquad
  \chi = \alpha + \beta_i c_i + \vartheta (|\bc|^2 - 3T). 
\end{equation}
Here, the Einstein summation convention is utilized. Thus, it holds that 
\begin{equation}
  \bv\cdot\nabla_{\bx}\mathcal{E}
  = \mathcal{M}\,\bv\cdot\nabla_{\bx}\chi
    + \chi\,\bv\cdot\nabla_{\bx}\mathcal{M}
  \triangleq  \mathcal{M}\bigl(\mathcal{T}_\rho + \mathcal{T}_u + \mathcal{T}_T\bigr),
  \label{eq:vgradH_decomp}
\end{equation}
with 
\begin{subequations}
\label{eq:modes}
\begin{align}
\label{eq:Trho_modes}
      \mathcal{T}_\rho
  &= \beta_k A_\ell\,\mathcal{A}_{k\ell}(\bc)
    + \vartheta A_k\,|\bc|^2 c_k
    + \mathcal{R}_{\rho},\\
      \label{eq:Tu_modes}
  \mathcal{T}_u
  &= S^{(u)}_{k\ell}\,\mathcal{A}_{k\ell}(\bc)
    + R^{(u)}_{k\ell}\,|\bc|^2\mathcal{A}_{k\ell}(\bc)
    + C^{(u,3)}_k\,|\bc|^2 c_k
    + \mathcal{R}_{u},\\
      \label{eq:TT_modes}
  \mathcal{T}_T
  &= S^{(T)}_{k\ell}\,\mathcal{A}_{k\ell}(\bc)
    + R^{(T)}_{k\ell}\,|\bc|^2\mathcal{A}_{k\ell}(\bc)
    + C^{(T,3)}_k\,|\bc|^2 c_k
    + C^{(T,5)}_k\,|\bc|^4 c_k
    + \mathcal{R}_{T},
\end{align}
\end{subequations}
and 
{\small
\begin{subequations}
    \begin{align}
  &S^{(u)}_{k\ell}
  = \left(\frac{\alpha}{T} - 3\vartheta\right)U_{k\ell}
    + \frac{D_{k\ell}}{T}
    - \frac{\delta u_\ell\,G_k}{T^2}
    + \frac{\beta_k\,u_{m}U_{m\ell}}{T},
  \qquad
  R^{(u)}_{k\ell} = \frac{\vartheta}{T}\,U_{k\ell},
  \label{eq:Su_Ru_def}\\
 & C^{(u,3)}_k
  = \frac{\vartheta}{T}\,u_{\ell}U_{\ell k}
    + \frac{1}{5T}\bigl(
        b_\ell U_{\ell k} + b_\ell U_{k\ell} + \beta_k U_{\ell\ell}
      \bigr),
  \label{eq:Cu3_def}\\
  &S^{(T)}_{k\ell}
  = -\frac{3b_k G_\ell}{2T}
     - \frac{\delta T\,U_{k\ell}}{T^2},
  \qquad
  R^{(T)}_{k\ell} = \frac{\beta_k G_\ell}{2T^2},
  \label{eq:ST_RT_def}\\
  &C^{(T,3)}_k
  = \frac{\alpha\,G_k}{2T^2}
     - \frac{3\vartheta\,G_k}{T}
     + \frac{\partial_k\delta T}{2T^2}
     - \frac{\delta T\,G_k}{T^3}
     + \frac{\beta_k\,u_{\ell}G_\ell}{2T^2},
  \qquad
  C^{(T,5)}_k = \frac{\vartheta\,G_k}{2T^2}.
  \label{eq:CT3_CT5_def}
\end{align}
\end{subequations}}
Here, only the terms that will contribute when calculating $\sigma^{(\nabla)}$ and $q^{(\nabla)}$ are listed, namely 
\begin{equation}
  \bigl\{\,
    \mathcal{A}_{k\ell}(\bc),\quad
    |\bc|^2\,\mathcal{A}_{k\ell}(\bc),\quad
    |\bc|^2\,c_k,\quad
    |\bc|^4\,c_k
  \,\bigr\},
  \label{eq:relevant_modes}
\end{equation}
while all the others are summarized in $\mathcal{R}_{s},s = \rho, u, T$. Moreover, due to the isotropy of the linearized collision operator $\mL$ \cite{cai2015approximation}, we can deduce that the integrals have the form below 
{\small
\begin{subequations}
\label{eq:def_coe}
    \begin{align}
  &\int_{\mathbb{R}^3}\mathcal{A}_{ij}(\bc)\,
    \mathcal{L}^{-1}\bigl[\mathcal{A}_{k\ell}(\bc)\,\mathcal{M}\bigr]
  \,\mathrm{d}\bv
  =: \frac{C_{2,0}}{\nu}
    \!\left(\delta_{ik}\delta_{j\ell}+\delta_{i\ell}\delta_{jk}
    -\tfrac{2}{3}\delta_{ij}\delta_{k\ell}\right),
  \label{eq:def_C20} \\
  &\int_{\mathbb{R}^3}\mathcal{A}_{ij}(\bc)\,
    \mathcal{L}^{-1}\bigl[|\bc|^2\mathcal{A}_{k\ell}(\bc)\,\mathcal{M}\bigr]
  \,\mathrm{d}\bv
  =: \frac{C_{2,1}}{\nu}\,\rho T^2
    \!\left(\delta_{ik}\delta_{j\ell}+\delta_{i\ell}\delta_{jk}
    -\tfrac{2}{3}\delta_{ij}\delta_{k\ell}\right),
  \label{eq:def_C21} \\
 & \int_{\mathbb{R}^3}\mathcal{B}_i(\bc)\,
    \mathcal{L}^{-1}\bigl[|\bc|^2 c_k\,\mathcal{M}\bigr]
  \,\mathrm{d}\bv
  =: -\frac{C_{1,1}}{\nu}\,\rho T^3\,\delta_{ik},
  \label{eq:def_C11} \\
 & \int_{\mathbb{R}^3}\mathcal{B}_i(\bc)\,
    \mathcal{L}^{-1}\bigl[|\bc|^4 c_k\,\mathcal{M}\bigr]
  \,\mathrm{d}\bv
  =: -\frac{C_{1,2}}{\nu}\,\rho T^4\,\delta_{ik},
  \label{eq:def_C12} \\
  &\int_{\mathbb{R}^3}\mathcal{A}_{ij}(\bc)\,
    \mathcal{L}^{-1}\bigl[|\bc|^2 c_k\,\mathcal{M}\bigr] = 0, \qquad \int_{\mathbb{R}^3}\mathcal{A}_{ij}(\bc)\,
    \mathcal{L}^{-1}\bigl[|\bc|^4 c_k\,\mathcal{M}\bigr] = 0, \\
  &   \int_{\mathbb{R}^3}\mathcal{B}_i(\bc)\,
    \mathcal{L}^{-1}\bigl[\mathcal{A}_{k\ell}(\bc)\,\mathcal{M}\bigr]
  \,\mathrm{d}\bv = 0, \qquad \int_{\mathbb{R}^3}\mathcal{B}_i(\bc)\,
    \mathcal{L}^{-1}\bigl[|\bc|^2\mathcal{A}_{k\ell}(\bc)\,\mathcal{M}\bigr]
  \,\mathrm{d}\bv = 0,
\end{align}
\end{subequations}}
where $\nu$ is the local collision frequency and $C_{2,0}$, $C_{2,1}$, $C_{1,1}$, $C_{1,2}$ are scalar coefficients determined by the linearized collision operator. Therefore, with the definition 
\begin{equation}
  \sigma^{(\nabla)}_{ij} 
  = \int_{\mathbb{R}^3} \mathcal{A}_{ij}(\bc)\, g_1^{(\nabla)} \,\mathrm{d}\bv,
  \qquad
  q^{(\nabla)}_i 
  = \int_{\mathbb{R}^3} \mathcal{B}_i(\bc)\, g_1^{(\nabla)} \,\mathrm{d}\bv,
  \label{eq:constitutive_def}
\end{equation}
it holds that 
{\small
\begin{subequations}
\label{eq:sigma_q_rho}
\begin{align}
  &\sigma^{(\rho)}_{ij}\triangleq \int_{\mathbb{R}^3}\mathcal{A}_{ij}(\bc)\,
    \mathcal{L}^{-1}\bigl[\mathcal{T}_{\rho}\,\mathcal{M}\bigr]
    \,\mathrm{d}\bv 
  = \beta_k A_\ell 
    \int_{\mathbb{R}^3}\mathcal{A}_{ij}(\bc)\,
    \mathcal{L}^{-1}\bigl[\mathcal{A}_{k\ell}(\bc)\,\mathcal{M}\bigr]
    \,\mathrm{d}\bv
  = \frac{C_{2,0}}{\nu}
    \bigl(\tf{\beta \otimes A}\bigr)_{ij},
  \label{eq:sigma_rho} \\
 &  q^{(\rho)}_i \triangleq \int_{\mathbb{R}^3}\mathcal{B}_{i}(\bc)\,
    \mathcal{L}^{-1}\bigl[\mathcal{T}_{\rho}\,\mathcal{M}\bigr]
    \,\mathrm{d}\bv 
  = \vartheta A_k
    \int_{\mathbb{R}^3}\mathcal{B}_i(\bc)\,
    \mathcal{L}^{-1}\bigl[|\bc|^2 c_k\,\mathcal{M}\bigr]
    \,\mathrm{d}\bv
  = -\frac{C_{1,1}}{\nu}\rho T^3\,\vartheta A_i.
  \label{eq:q_rho} 
  \end{align}
\end{subequations}}
Similarly, we can obtain that 
{\small
\begin{subequations}
\label{eq:sigma_q_u}
\begin{align}
  & \sigma^{(u)}_{ij}
  \triangleq \int_{\mathbb{R}^3}\mathcal{A}_{ij}(\bc)\,
    \mathcal{L}^{-1}\bigl[\mathcal{T}_{u}\,\mathcal{M}\bigr]
    \,\mathrm{d}\bv 
  = \frac{C_{2,0}}{\nu}\left(\tf{S^{(u)}}\right)_{ij}
  + \frac{C_{2,1}}{\nu}\rho T^2\left(\tf{R^{(u)}}\right)_{ij},   
  \label{eq:sigma_u_result}\\
 & q^{(u)}_i \triangleq \int_{\mathbb{R}^3}\mathcal{B}_{i}(\bc)\,
    \mathcal{L}^{-1}\bigl[\mathcal{T}_{u}\,\mathcal{M}\bigr]
    \,\mathrm{d}\bv =  -\frac{C_{1,1}}{\nu}\rho T^3\,C^{(u,3)}_i,
  \label{eq:q_u_result}
  \end{align}
\end{subequations}}
and 
{\small
\begin{subequations}
\label{eq:sigma_q_T}
\begin{align}
  &  \sigma^{(T)}_{ij}
  \triangleq \int_{\mathbb{R}^3}\mathcal{A}_{ij}(\bc)\,
    \mathcal{L}^{-1}\bigl[\mathcal{T}_{T}\,\mathcal{M}\bigr]
    \,\mathrm{d}\bv = \frac{C_{2,0}}{\nu}\left(\tf{S^{(T)}}\right)_{ij}
  + \frac{C_{2,1}}{\nu}\rho T^2\left(\tf{R^{(T)}}\right)_{ij},
  \label{eq:sigma_T_result} \\
  & q^{(T)}_i
  \triangleq \int_{\mathbb{R}^3}\mathcal{B}_{i}(\bc)\,
    \mathcal{L}^{-1}\bigl[\mathcal{T}_{T}\,\mathcal{M}\bigr]
    \,\mathrm{d}\bv
  = -\frac{C_{1,1}}{\nu}\rho T^3\,C^{(T,3)}_i
    -\frac{C_{1,2}}{\nu}\rho T^4\,C^{(T,5)}_i.
  \label{eq:q_T_result}
  \end{align}
\end{subequations}}
Collecting \eqref{eq:sigma_q_rho}, \eqref{eq:sigma_q_u}, and \eqref{eq:sigma_q_T}, we finally achieve that 
\begin{align}
  \sigma^{(\nabla)}
  = \frac{C_{2,0}}{\nu}\,\mathbf{S}
   + \frac{C_{2,1}}{\nu}\,\mathbf{R},
    \qquad 
  q^{(\nabla)}
  = -\frac{1}{2\nu}\bigl(C_{1,1}\,\boldsymbol{P} + C_{1,2}\,\boldsymbol{Q}\bigr),
  \end{align}
where $\mathbf{S}$, $\mathbf{R}$, $\boldsymbol{P}$, $\boldsymbol{Q}$ are defined in \eqref{eq:eq_s} and \eqref{eq:eq_q}.

\subsection{Calculating the coefficients \texorpdfstring{$C_{i,j}$}{C{i,j}}}
\label{app:bgk_coefficients}
We first present the integrals of the Maxwellian as below, which are utilized to calculate $C_{i,j}$.
\begin{gather}
  \label{eq:gauss}
   \int_{\bbR^3} c_ic_j \mathcal{M}\,\mathrm{d}\bv 
  = \rho T\delta_{ij},
  \qquad 
  \int_{\bbR^3} c_ic_jc_kc_\ell \mathcal{M}\,\mathrm{d}\bv 
  = \rho T^2(\delta_{ij}\delta_{k\ell}
    +\delta_{ik}\delta_{j\ell}+\delta_{i\ell}\delta_{jk}), \\
   \int_{\bbR^3} |\bc|^2c_ic_j \mathcal{M}\,\mathrm{d}\bv 
  = 5\rho T^2\delta_{ij},\quad 
  \int_{\bbR^3} |\bc|^4c_ic_j \mathcal{M}\,\mathrm{d}\bv 
  = 35\rho T^3\delta_{ij},
  \quad
  \int_{\bbR^3} |\bc|^6c_ic_j \mathcal{M}\,\mathrm{d}\bv 
  = 315\rho T^4\delta_{ij}.
\end{gather}
To obtain  the coefficients $C_{i,j}$ in \eqref{eq:def_coe}, the technique of utilizing the BGK operator $\mL_{\rm BGK}$ instead of the linearized collision operator $\mL$ is applied. The pseudoinverse of $\mL_{\rm BGK}$ has the form as 
\begin{equation}
  \mathcal{L}_{\mathrm{BGK}}^{-1}[g] = -\frac{1}{\nu}(I-\mathcal{P})g.
  \label{eq:BGK_inv_1}
\end{equation}
Then, it is straightforward to verify 
\begin{equation}
    \label{eq:orth_A}
    \int_{\bbR^3} \phi(\bv) \mathcal{A}_{k\ell}\mathcal{M} \dd \bv = 0, \qquad \phi(\bv) = (1, \bv, |\bv|^2)^T,
\end{equation}
so that 
\begin{equation}
\label{eq:c20_L}
    (I-\mathcal{P})[\mathcal{A}_{k\ell}\mathcal{M}] = \mathcal{A}_{k\ell}\mathcal{M},\qquad \mathcal{L}_{\mathrm{BGK}}^{-1}[\mathcal{A}_{k\ell}\mathcal{M}] = -\frac{1}{\nu}\mathcal{A}_{k\ell}\mathcal{M}.
\end{equation}
Substituting \eqref{eq:c20_L} into \eqref{eq:def_C20}, with \eqref{eq:gauss} it yields
\begin{align}
  \int_{\bbR^3}\mathcal{A}_{ij}\,\mathcal{L}_{\mathrm{BGK}}^{-1}[\mathcal{A}_{k\ell}\mathcal{M}]
  \,\mathrm{d}\bv
  &= -\frac{1}{\nu}\int\mathcal{A}_{ij}\mathcal{A}_{k\ell}
    \mathcal{M}\,\mathrm{d}\bv 
  = -\frac{\rho T^2}{\nu}
    \!\left(\delta_{ik}\delta_{j\ell}+\delta_{i\ell}\delta_{jk}
    -\frac{2}{3}\delta_{ij}\delta_{k\ell}\right).
    \label{eq:C20}
\end{align}
Comparing \eqref{eq:C20} with \eqref{eq:def_C20}, it gives
\begin{equation}
    \label{eq:C_20}
  C_{2,0}=-\rho T^2.    
\end{equation}
Similarly, it holds for $|\bc|^2\mathcal{A}_{k\ell}\mathcal{M}$ that 
\begin{equation}
    \label{eq:orth_bA}
      \int_{\bbR^3} \phi(\bv) |\bc|^2\mathcal{A}_{k\ell}\mathcal{M} \dd \bv = 0,
\end{equation}
with \eqref{eq:BGK_inv_1} and \eqref{eq:gauss}, we can derive that 
\begin{align}
  \int_{\bbR^3} \mathcal{A}_{ij}\,
  \mathcal{L}_{\mathrm{BGK}}^{-1}\bigl[|\bc|^2\mathcal{A}_{k\ell}\mathcal{M}\bigr]
  \,\mathrm{d}\bv
    = -\frac{7\rho T^3}{\nu}
  \left(\delta_{ik}\delta_{j\ell}+\delta_{i\ell}\delta_{jk}
  -\frac{2}{3}\delta_{ij}\delta_{k\ell}\right).
  \label{eq:C21}
\end{align}
Comparing \eqref{eq:C21} with \eqref{eq:def_C21}, it gives 
\begin{equation}
    \label{eq:C_21}
   C_{2,1} = -7T.
\end{equation}
Moreover, for $|\bc|^2c_k\mathcal{M}$, and $|\bc|^4c_k\mathcal{M}$, we have that 
\begin{equation}
\label{eq:bcM}
\int_{\bbR^3} c_i|\bc|^2c_k\mathcal{M}\,\mathrm{d}\bv = 5\rho T^2\delta_{ik}, \qquad  \int_{\bbR^3} c_i\cdot|\bc|^4c_k\mathcal{M}\,\mathrm{d}\bv 
  = 35\rho T^3\delta_{ik}.
\end{equation}
Then, it holds that 
\begin{equation}
    \mathcal{P}[|\bc|^2c_k\mathcal{M}] = 5T c_k\mathcal{M}, \qquad \mathcal{P}[|\bc|^4c_k\mathcal{M}] = 35T^2\,c_k\mathcal{M}.
    \label{eq:C11_L}
\end{equation}
Using \eqref{eq:BGK_inv_1},\eqref{eq:gauss}, \eqref{eq:bcM} and \eqref{eq:C11_L}, we obtain
\begin{equation}
    \int_{\bbR^3} \mathcal{B}_i\, \mathcal{L}_{\mathrm{BGK}}^{-1}[|\bc|^2c_k\mathcal{M}] \,\mathrm{d}\bv
    = -\frac{5\rho T^3}{\nu}\delta_{ik}, \qquad  \int_{\bbR^3} \mathcal{B}_i\,\mathcal{L}_{\mathrm{BGK}}^{-1}[|\bc|^4c_k\mathcal{M}]
  \,\mathrm{d}\bv
   = -\frac{70\rho T^4}{\nu}\,\delta_{ik}.
    \label{eq:C11}
\end{equation}
With \eqref{eq:C11} and \eqref{eq:def_coe}, it gives
\begin{equation}
    \label{eq:C11_C12}
    C_{1,1} = 5, \qquad C_{1,2} = 70.
\end{equation}

In summary, with the BGK approximation, the coefficients $C_{i,j}$ take explicit values as 
\begin{equation}
  C_{2,0} = -\rho T^2,\quad
  C_{2,1} = -7T,\quad
  C_{1,1} = 5,\quad
  C_{1,2} = 70.
  \label{eq:coeff_values}
\end{equation}
\bibliographystyle{plain}
\bibliography{article}

\end{document}